\documentclass[a4paper, 10pt]{article}

\usepackage[a4paper,left=3.3cm,right=3.3cm,top=2.8cm,bottom=2.8cm]{geometry}

\usepackage{pgfplots}
\usepackage{pgfplotstable}
\usepgfplotslibrary{groupplots}
\usepgfplotslibrary{fillbetween}
\usepackage{mathtools}
\usepackage{multicol}
\usepackage{multirow}
\usepackage{comment}
\usepackage{booktabs}
\pgfplotsset{compat=1.5}
\usepackage{amssymb}
\usepackage{url}
\usepackage{bm}
\usepackage{algorithm}
\usepackage{algpseudocode}

\usepackage{pdfsync}
\usepackage{float}
\usepackage{tabularx}
\usepackage{enumerate}
\usepackage{array}
\usepackage{xspace}
\usepackage{tikz}
\usepackage{tikzsymbols}
\usetikzlibrary{calc,trees,positioning,arrows,chains,shapes.geometric,%
    decorations.pathreplacing,decorations.pathmorphing,shapes,%
    matrix,shapes.symbols, decorations.markings, patterns,fit}
\usepackage{color, colortbl}

\usepackage{subcaption}
\usepackage{authblk}
\usepackage{units}
\usepackage[scientific-notation=true]{siunitx}

\usepackage{hyperref}
\usepackage[draft,inline,marginclue]{fixme}
\FXRegisterAuthor{mt}{amt}{MT}

\usepackage{cleveref}

\newcommand{\mupar}{\ensuremath{\boldsymbol{\mu}}}
\newcommand{\etapar}{\ensuremath{\boldsymbol{\eta}}}

\newcommand{\X}{\mathbf{X}}
\newcommand{\R}{\mathbb{R}}
\newcommand{\x}{\mathbf{x}}
\newcommand{\y}{\mathbf{y}}

\newcommand{\RA}[1]{{\color{black}#1}}
\newcommand{\RB}[1]{{\color{black}#1}}

\definecolor{Gray}{gray}{0.9}

\makeatletter
\def\algbackskip{\hskip-\ALG@thistlm}
\makeatother

\begin{document}

\title{A supervised learning approach involving active subspaces for an
efficient genetic algorithm in high-dimensional optimization problems} 

\author[]{Nicola~Demo\footnote{nicola.demo@sissa.it}}
\author[]{Marco~Tezzele\footnote{marco.tezzele@sissa.it}}
\author[]{Gianluigi~Rozza\footnote{gianluigi.rozza@sissa.it}}

\affil{SISSA, Mathematics Area, mathLab, via Bonomea 265, I-34136
  Trieste, Italy}

\maketitle

\begin{abstract}

In this work, we present an extension of the genetic algorithm (GA) which
exploits the \RA{supervised learning technique called active subspaces (AS)} to evolve the individuals on a
lower dimensional space. In many cases, GA requires in fact more function
evaluations than others optimization method to converge to the
\RA{global} optimum. Thus, complex and high-dimensional functions may
	result \RA{extremely demanding (from computational viewpoint) to optimize} with the standard algorithm. To
address this issue, we propose to linearly map the input parameter
space of the original function onto its AS before the evolution,
performing the \textit{mutation} and \textit{mate} processes in a
lower dimensional space. In this contribution, we describe the novel
method called ASGA, presenting differences and similarities with the
standard GA method. We test the proposed method over $n$-dimensional
benchmark functions --- {\it Rosenbrock}, {\it Ackley}, {\it
  Bohachevsky}, {\it Rastrigin}, {\it Schaffer N.~7}, and {\it
  Zakharov} --- and finally we apply it to an aeronautical shape
optimization problem.

\end{abstract}

\tableofcontents

\section{Introduction}
\label{sec:intro}
Genetic algorithm (GA) is a well-known and widespread methodology, mainly
adopted in optimization problems~\cite{holland1992adaptation,kumar2010genetic}. It emulates the
evolutive process of \RA{natural selection} by following an iterative process where
the individuals are selected by a given objective function and subsequently
they mutate and reproduce~\cite{emmeche1996garden, ali2005numerical,
chelouah2000continuous, laguna2005experimental}. This gradient-free technique
is particular effective when the objective function contains many local minima:
thanks to the stochastic component, GA explores the domain without being blocked into
 local minima.  The main disadvantages of such algorithm is the (relative)
high number of required evaluations of the objective function during the
evolution \RA{to explore the input space~\cite{nesterov2017random}}, that makes in
several industrial and engineering contexts this method unfeasible for the
global computational cost.

In this work, we propose a novel extension of standard GA, exploiting the
emerging active \RA{subspaces} (AS) property~\cite{constantine2015active,
constantine2014active} for the dimensionality reduction. \RA{AS is a
supervised learning technique which} allows to
approximate a scalar function with a lower dimensional one, whose parameters
are a linear combination of the original inputs. 
AS has been successfully employed in naval engineering
applications~\cite{tezzele2018dimension, tezzele2018ecmi, tezzele2018model,
tezzele2019marine}, coupled with reduced order methods such as POD-Galerkin in
biomedical applications~\cite{tezzele2018combined,morhandbook2019}, POD with
interpolation~\cite{demo2019cras} in structural and CFD analysis, and Dynamic
Mode Decomposition~\cite{tezzele2019mortech} in CFD contexts. Other
applications \RA{include} aerodynamic shape
optimization~\cite{lukaczyk2014active}, artificial neural networks to
reduce the number of neurons~\cite{cui2019active}, non-linear
structural analysis~\cite{guo2018reduced}, and AS for multivariate
vector-valued model functions~\cite{zahm2020gradient}. Several
non-linear AS extension have been proposed recently. We mention Active
Manifold~\cite{bridges2019active}, Kernel-based Active Subspaces
(KAS)~\cite{romor2020kas} which exploits the random Fourier features
to map the inputs in a higher dimensional space. We also mention the application of 
artificial neural networks for non-linear reduction in parameter spaces by
learning isosurfaces~\cite{zhang2019learning}.  Despite these new non-linear
extensions of AS, in this work we exploit the \RA{classical linear} version because
of the possibility to map points in the reduced space onto the original
parameter space.

The main idea of the proposed algorithm is to force the individuals of the
population to evolve along the AS, which has a lower dimension, avoiding evolution
along the meaningless directions. Further,
the high number of function evaluations that characterize the GA is exploited
within this new approach for the construction (and refinement) of the AS,
making these techniques --- GA produces a large dataset of input-output pairs,
whereas AS needs large datasets for an accurate subspace identification ---
particularly suited together. This new method has the potential to
improve existing optimization pipeline involving both input and model
order reduction.

\RA{A similar approach has been proposed in~\cite{choromanski2019complexity}, where an active subspace
is constructed in order to obtain an efficient and adaptive sampling
strategy in a evolution strategy framework. This approach shares with
the one we are proposing the idea of efficiently exploring the input
space by constructing a subspace based on the collected data. In
contrast with our approach in~\cite{choromanski2019complexity} the subspace construction is done
with a singular value decomposition based method, and the
optimization technique is completely different, even if evolution
strategy methods and genetic algorithm present some analogies. To the best of the authors
knowledge, the current contribution presents a novel approach, not yet
explored in the literature. For similar approach, we cite also}
random subspace embeddings for unconstrained
global optimization of functions with low effective dimensionality that can be
found in~\cite{wang2016bayesian,cartis2020dimensionality}, while for
evolutionary methods and derivative-free optimization we
mention~\cite{sanyang2016remeda, qian2016derivative}, respectively. For a
survey on linear dimensionality reduction in the context of optimization
programs over matrix manifolds we mention~\cite{cunningham2015linear}.

The outline of this work is the following: the proposed method is
described in \cref{sec:gas}, while \cref{sec:genetic} and
\cref{sec:active} are devoted to recall the general family of
genetic algorithms and the active subspaces technique,
respectively. \cref{sec:results} presents the numerical results
obtained applying the proposed extension to some popular benchmark
functions for optimization problems, then to a typical engineering
problem where the shape of a NACA airfoil is morphed to maximize the
lift-to-drag coefficient. Finally \cref{sec:conclusions} 
summarizes the benefits of the method and proposes some extensions for
future developments.

\section{Genetic Algorithms}
\label{sec:genetic}

In this work we propose an extension of the standard genetic
algorithm (GA). We start recalling the general method in order to easily
let the reader understand the differences. We define GA as the family of
computational methods that are inspired by Darwin's theory of evolution. The
basic idea is to generate a population of individuals with random genes, and
make them evolve through mutations and crossovers, mimicking the evolution of
living beings. Iterating this process by selecting at each step the best-fit
individuals results in the optimization --- according to a specific objective
function --- of the original population. As such this method can be
easily adopted as a global optimization algorithm.

Initially proposed by Holland in~\cite{holland1973genetic}, GA has
had several modifications during the
years (see for example~\cite{kumar2010genetic,el2006hybrid,sivaraj2011review,elsayed2014new,drezner2003new}),
but it keeps its fundamental steps: \textit{selection}, 
\textit{mutation} and \textit{mate}.

\begin{figure}[ht]
\begin{subfigure}{.49\textwidth}
  \centering
  \includegraphics[width=.9\linewidth, trim=20 0 0 0, clip]{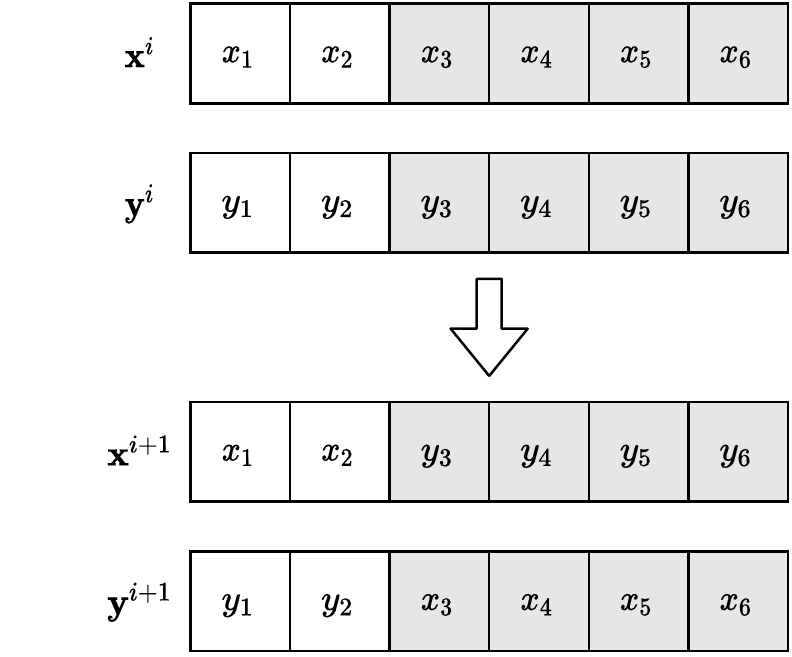}
  \caption{Mate}
  \label{fig:crossover}
\end{subfigure}
\begin{subfigure}{.49\textwidth}
  \centering
  \includegraphics[width=.9\linewidth, trim=20 0 0 0, clip]{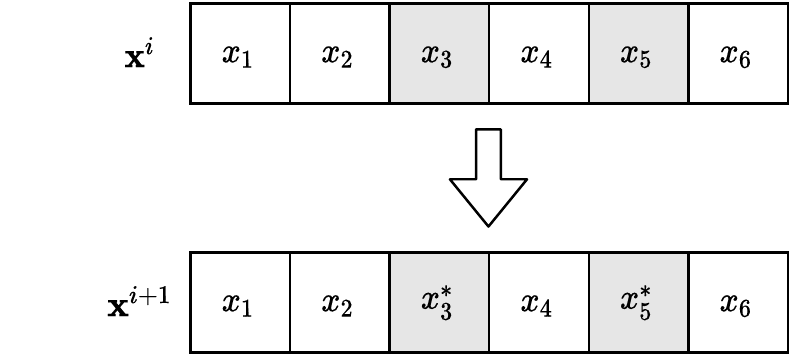}
  \caption{Mutation}
  \label{fig:mutation}
\end{subfigure}
	\caption{Graphical example of mate and mutation \RB{where
            $\mathbf{x}^i$ and $\mathbf{y}^i$ indicates two generic
            individuals of the $i$-th generation}.}
\label{fig:fig}
\end{figure}

Let us define formally the individuals: a population
composed by $N$ individuals $\x_i \in \R^P$ with $P$ genes is defined
as $\X = \{\x_1, \dotsc, \x_N\}$. We express the fitnesses of such
individuals with the scalar function $f: \R^P \to \R$. The first generation
$\X^1$ is randomly created --- with possible constraints --- and the
fitness is evaluated for all the individuals: $\y_i = f(\x_i)$ for $i
= 1, \dotsc, N$. Then the following iterative process starts:
\begin{description}
	\item[Selection:] The best individuals of the previous generation
		$\X^i$ are chosen accordingly to their fitnesses to breed the
		new generation. For the selection, several strategies can be
		adopted depending on the problem and on the \RA{cardinality
		of the population $N$}. 
	\item[Mate:] Finally, the selected individuals are grouped into pairs
		and, according to a mate probability, they combine their genes
		to create new individuals. The process, also called {\it
		crossover} emulate the species reproduction. These individuals
		form the new generation $\X^{i+1}$. An example of a crossover
		method is sketched in \cref{fig:crossover}.
	\item[Mutation:] The individuals evolve by changing some of their
		genes. The mutation of an individual is usually controlled by a
		mutation probability. In \cref{fig:mutation} we show an
		illustrative example where two genes have randomly mutated.
\end{description}

After the mutation step, the fitness of the new individuals is computed
and the algorithm restarts with the selection of the best-fit
individuals. In this way, the population evolves, generation after
generation, towards the optimal individual, \RA{avoiding getting blocked} in a
local minima thanks to the stochastic component introduced by mutation and
crossover. Thus, this method is very effective for global optimization
where the objective function is potentially non-linear, while standard
gradient-based methods can converge to local minima. However GA usually
requires an high number of evaluations to perform the optimization, making this
procedure very expensive in \RA{the} case of computational costly objective functions.

\section{Active subspaces for minimization on a lower dimensional
  parameter space}
\label{sec:active}
Active Subspaces (AS)~\cite{constantine2015active}
method is a dimensionality reduction approach for parameter space
studies\RA{, which falls under the category of supervised learning
  techniques}. AS tries to reduce the input dimension of a scalar
function $f(\mupar) : \Omega \subset \mathbb{R}^k \to \mathbb{R}$  by defining a linear
transformation $\mupar_M = \mathbf{A} \mupar$. This method requires
the evaluation of the gradients of $f$ 
since $\mathbf{A}$ depends on the second moment matrix $\mathbf{C}$ of
the target function’s gradient, also called uncentered covariance matrix of the
gradients of $f$. This matrix is defined as follows
\begin{equation}
\label{eq:cov_matrix}
\mathbf{C} = \mathbb{E}\, [\nabla_{\mupar} f \, \nabla_{\mupar} f
^T] =\int (\nabla_{\mupar} f) ( \nabla_{\mupar} f )^T
\rho \, d \mupar,
\end{equation}
where with the symbol $\mathbb{E} [\cdot]$ we denote the expected
value, $\nabla_{\mupar} f \equiv \nabla f({\mupar}) \in \mathbb{R}^k$,
and $\rho : \mathbb{R}^k \to \mathbb{R}^+$ is a probability density
function which represents the uncertainty in the input parameters.  In
practice the matrix $\mathbf{C}$ is constructed with a Monte Carlo
procedure, and the gradients if not provided can be approximated with
different techniques such as local linear models, global models,
finite difference, or Gaussian process for example \RB{--- for a
  comparison of the methods and corresponding errors the reader can refer
  to~\cite{constantine2015active, constantine2014active, constantine2015computing}}. The
uncentered covariance matrix can be decomposed~as
\begin{equation}
\mathbf{C} = \mathbf{W} \Lambda \mathbf{W}^T,
\end{equation}
where $\mathbf{W}$ stands for the orthogonal matrix containing the
eigenvectors, and $\Lambda$ for the eigenvalues matrix arranged in
descending order. The \RA{spectral gap~\cite{constantine2015active}} is used to bound the error on the
numerical approximation with Monte Carlo. We can decompose these
two matrices as
\begin{equation}
\Lambda =   \begin{bmatrix} \Lambda_1 & \\
& \Lambda_2\end{bmatrix},
\qquad
\mathbf{W} = \left [ \mathbf{W}_1 \quad \mathbf{W}_2 \right ],
\qquad
\mathbf{W}_1 \in \mathbb{R}^{k\times M},
\end{equation}
where $M < k$ is the dimension of the active subspace. $M$ should be
chosen looking at the energy decay (the tail in the ordered
eigenvalue sum) as in POD, or it can be prescribed a priori for the
specific task. We can exploit this decomposition to map the input
parameters onto a reduced space.

We define the active subspace of dimension $M$ as
the principal eigenspace corresponding to the eigenvalues prior to the
major spectral gap. We also call the active variable $\mupar_M$ and
the inactive variable $\etapar$. They are defined as $\mupar_M =
\mathbf{W}_1^T\mupar \in \mathbb{R}^M$, and $\etapar =
\mathbf{W}_2^T\mupar \in \mathbb{R}^{k-M}$.

In this work we address the constrained global optimization problem of
a real-valued continuous function, in the context of genetic
algorithms, defined as
\begin{equation}
\min_{\mupar \in \Omega \subset \mathbb{R}^k} f(\mupar).
\end{equation}
To fight the curse of dimensionality and speed up the convergence we
exploit the active subspaces property of the target function to select
the best individuals in the reduced parameter space, mutate and mate
them, and successively to map them in the full parameter space. This
translates in the following optimization problem for each generation
of individuals: 
\begin{equation}
  \min_{\substack{\mupar_M \in \mathcal{P} \subset \mathbb{R}^M  \\
      \mupar \in \Omega}} \; g(\mupar_M = \mathbf{W}_1^T\mupar), 
\end{equation}
where $\mathcal{P}$ is the polytope in $\mathbb{R}^M$ --- we assume
the ranges of the parameters to be intervals --- defined by the
AS as $\mathcal{P} := \{ \mupar_M =  \mathbf{W}_1^T\mupar \; | \; \mupar \in
\Omega \}$. We remark that there are many choices for the profile $g$. In
this work we consider the following profile:
\begin{gather}
  g(\mathbf{y}) := f(\mathbf{x_y}) \qquad \forall \mathbf{y} \in \mathcal{Y}, \\
  \mathcal{Y} := \{ \mathbf{y} \in \mathcal{P} \; | \; \exists \, \mathbf{x_y} \in \Omega \,
  \text{ s.t. } \mathbf{y} =  \mathbf{W}_1^T \mathbf{x_y} \}.
\end{gather}

\begin{figure}[ht]
\centering
\includegraphics[width=1.\linewidth]{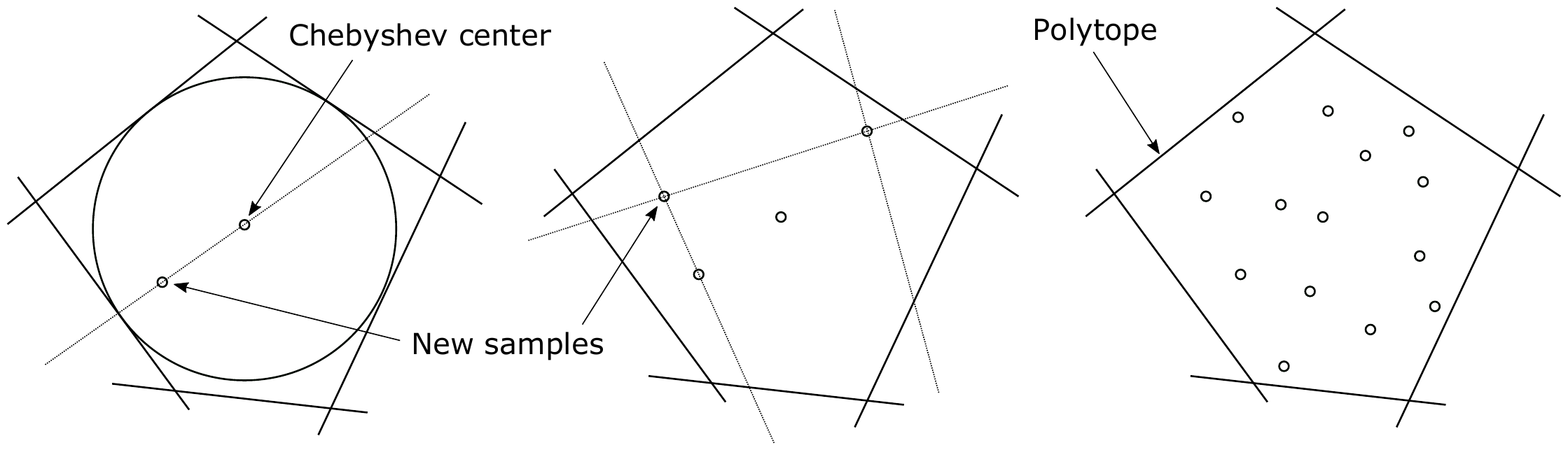}
\caption{\RA{Graphical representation of the inactive variable sampling
  strategy. We emphasize the Chebyshev center, the selection of the
  next sample using the \emph{hit and run} method, and the polytope defined
  by \cref{eq:polytope_def}.\label{fig:hit_and_run}}}
\end{figure}

We \RA{emphasize} that the projection map onto the active subspace is a
surjective map because $\mathbf{W}_1^T$ is defined as a linear
projection onto a subspace, hence it is surjective by definition. So the
back-mapping from the active subspace onto $\Omega$ is not
trivial. Given a point $\mupar_M^*$ in the active subspace we can find
$B$ points in the original parameter space which are mapped onto
$\mupar_M^*$ by $\mathbf{W}_1^T$. \RA{Recalling the decomposition
  above we have that
  \begin{equation}
    \mupar = \mathbf{W}_1 \mathbf{W}_1^T \mupar +
    \mathbf{W}_2 \mathbf{W}_2^T \mupar = \mathbf{W}_1 \mupar_M +
    \mathbf{W}_2 \etapar \qquad \forall \mupar \in \Omega,
  \end{equation}
with the additional constraint coming from the rescaling of the input
parameters needed to apply AS: $-\mathbf{1} \leq \mupar \leq
\mathbf{1}$, where $\mathbf{1}$ denotes the vector in $\R^k$ with all
elements equal to $1$. We exploit this to sample the inactive variable $\etapar$ so that
\begin{equation}
-\mathbf{1} \leq \mathbf{W}_1 \mupar_M^* + \mathbf{W}_2 \etapar \leq \mathbf{1},
\end{equation}
or equivalently
\begin{equation}
  \label{eq:polytope_def}
 \mathbf{W}_2 \etapar \leq \mathbf{1} - \mathbf{W}_1 \mupar_M^* \qquad  - \mathbf{W}_2 \etapar
 \leq \mathbf{1} + \mathbf{W}_1 \mupar_M^*.
\end{equation}
}
These inequalities define a polytope in $\mathbb{R}^{k-M}$ from which
we want to uniformly sample $B$ points. \RA{The inactive variables are
  sampled from the conditional distribution $p(\etapar | \mupar_M^*)$,
  and we show how to perform it for the uniform distribution. For a
  more general distribution one should use Hamiltonian Monte Carlo.} In particular we start
with a simple rejection sampling scheme, which finds a bounding
hyperbox for the polytope, draws points uniformly from it, and rejects
points outside the polytope. If this method does not return enough
samples, we try a \RA{\emph{hit and run}
method~\cite{smith1996hit,belisle1993hit,lovasz2006hit} for sampling from the
polytope. This method, starting from the center of the largest
hypersphere within the polytope, selects a random direction and
identifies the longest segment lying inside the polytope. A new
sample is randomly drawn along this segment. The procedure continues
starting from the last sample until enough samples are found.
If also that does not work, we use $B$ copies of a feasible
point computed as the Chebyshev center~\cite{botkin1994algorithm} of
the polytope. In \cref{fig:hit_and_run} we depicted these strategies,
at different stages of the sampling.}

\section{The proposed ASGA optimization algorithm} 
\label{sec:gas}
\begin{figure}[ht]
\centering
\includegraphics[width=.7\linewidth]{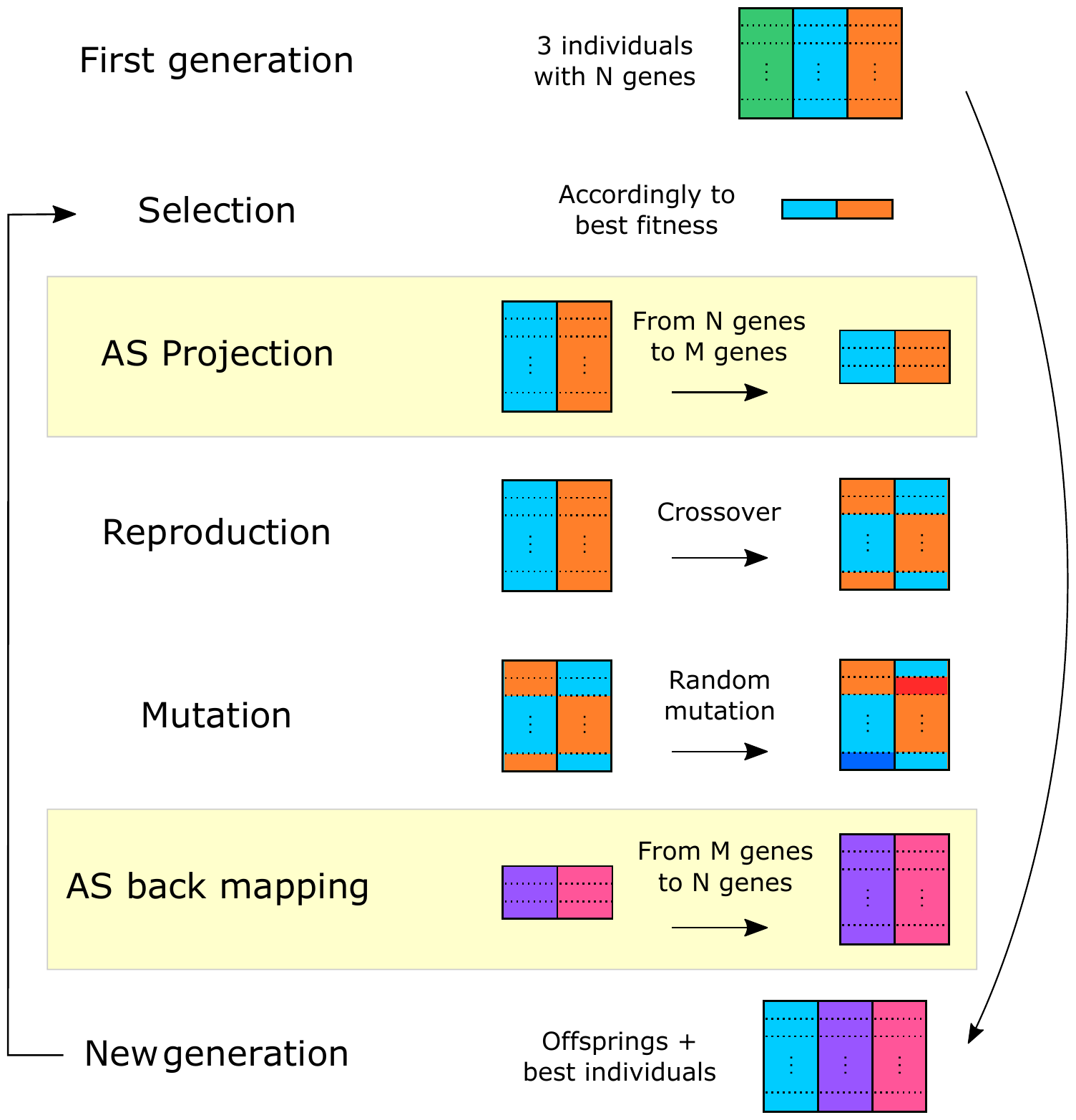}
\caption{\RA{Active subspaces-based genetic algorithm scheme. The main
    step of the classical GA are depicted from top to bottom. The
    yellow boxes represent projections onto and from lower dimension
    active subspace, which are specific to ASGA.\label{fig:asga}}} 
\end{figure}

In this section we are going to describe the proposed active subspaces extension of the
standard GA, named ASGA. Before
starting, we \RA{emphasize} that in what follows, we will maintain the
selection, mutation and mate procedures --- presented in
\cref{sec:genetic} --- as general as \RA{possible}, without going into
technical details, given the large variety of different options for these
steps. In fact the proposed extension is independent on the chosen
evolution strategies, and we only perform them in a lower dimension
exploiting AS. 
In \cref{algorithm} we summarize the standard approach, while
in \cref{algorithm1} we highlight the differences introduced
by ASGA. \RA{We also present an illustration for both the methods in
  \cref{fig:asga}, where the yellow boxes indicate the main steps
  peculiar to ASGA.} In both cases, the
first step is the generation of the random individuals composing the
initial population, and the sequential evaluation of all of them. For ASGA
these individuals and their fitness are stored into two
additional sets, $\mathbf{X}^{\text{AS}}$ for the individuals, and $\mathbf{y}^{\text{AS}}$
for the fitness. We will exploit them as input-output pair for
the construction of the AS. After the selection of the best-fit individuals,
the active subspace of dimension $M$ is built and the selected offspring is
projected onto it. The low-dimensional individuals mate and mutate
in the active subspace. Thanks to the reduced dimension and to the
fact that we retain only the most important dimensions, these operations
are much more efficient. Thus, even if the AS of dimension $M$ does
not provide an accurate approximation of the original full-dimensional space, the
active dimensions will provide preferential directions for the evolution,
making the iterative process smarter and faster.

After the evolution, the low-dimensional offspring is mapped back to the
original space. In \cref{sec:active} we describe how for any point in
the active subspace we can find several points in the original space
which are mapped onto it. So we select, for any individual in the
offspring, $B$ full-dimensional points which correspond to the
individual in the active subspace. We \RA{emphasize} that to
preserve the same dimensionality of the offspring between the original GA and
the AS extension, in the proposed algorithm we select the $\nicefrac{N}{B}$ best
individuals, instead of selecting $N$. In this way, after the back-mapping, the 
offspring has dimension $N$ in both versions. The number of back-mapped
points $B$, and the active subspace dimension $M$ --- that can be a
fixed parameter or dynamically selected from the spectral gap of the
covariance matrix $\mathbf{C}$ --- represent the new (hyper-)parameters
of the proposed method.

Finally the fitness of the new individuals, now in the
full-dimensional space, are evaluated. To make the AS more precise
during the iterations, the evaluated individuals 
and their fitness are added to $\mathbf{X}^{\text{AS}}$ and $\mathbf{y}^{\text{AS}}$.
The process restarts from the selection of the offspring from the new
generation, continuing as described above until the stopping criteria
are met.

\RB{We stress on the fact that the structure of the algorithm is
        similar to the
	original GA approach, with the difference that the gradients at the
	sample points are approximated in order to identify the dimensions with highest
	variance. Even if such information about the function gradient is used,
	the ASGA method is different from gradient-based methods: numerically
	computing the gradient with a good accuracy at a specific point  ---
	that is the fundamental step of gradient-based methods to move on the
	solution manifold --- is a very expensive procedure, especially in a
	high-dimensional space. In ASGA we avoid such computation, exploiting
	instead the already collected function evaluations. Further,
	gradient-based techniques converges (relatively) fast to optimum, but
	they get blocked into local minima, contrarily to the ASGA approach.}
It is important to remark that, for each generation, the AS is rebuilt
from scratch, \RA{losing} efficiency but gaining more precision due to the
growing number of elements in the two sets $\mathbf{X}^{\text{AS}}$ and
$\mathbf{y}^{\text{AS}}$. \RA{We also remark the samples are
generated with a uniform distribution only at the first
generation. After that, due to the ASGA steps the distribution
changes in a way which can not be known a priori. For the computation
of the expectation operator in \cref{eq:cov_matrix}
in this work we assume a uniform distribution. Even if this may
introduce an unknown error, the numerical results achieved 
by ASGA seem to support such choice. Of course the
numerical estimates present in the literature for the uniform
distribution do not apply in such case. This method can be viewed as an
active learning procedure in a Bayesian integration context, where the maximized
acquisition function is heuristic and given by the
application of AS and GA steps. Another interpretation is that we are
enriching the local informations near the current minimum to feed the
AS algorithm, so it can be viewed as a weighted AS.}

\newcommand{\red}[1]{\textcolor{black}{#1}}

\begin{minipage}{0.44\textwidth}
\begin{algorithm}[H]
    \scriptsize
    \caption{Standard GA.}\label{algorithm}
    \hspace*{\algorithmicindent} \textbf{Input}:\\
    \hspace*{\algorithmicindent}\hspace*{10pt} initial population size $N_0$\\
    \hspace*{\algorithmicindent}\hspace*{10pt} population size $N$\\
    \\
    \\
    \hspace*{\algorithmicindent}\hspace*{10pt} selection routine \textsc{select}\\
    \hspace*{\algorithmicindent}\hspace*{10pt} mutation routine \textsc{mutate}\\
    \hspace*{\algorithmicindent}\hspace*{10pt} mate routine \textsc{mate}\\
    \hspace*{\algorithmicindent}\hspace*{10pt} objective function \textsc{fobj}\\
    \hspace*{\algorithmicindent}\hspace*{10pt} stop criteria\\
    \hspace*{\algorithmicindent} \textbf{Output}:\\ 
    \hspace*{\algorithmicindent}\hspace*{10pt} final population $\mathbf{X}^{\text{end}}$\\
    \begin{algorithmic}[1]
        \Procedure{GeneticAlgorithm}{}
        \State $g \gets$ 0 
        \State $\mathbf{X}^g \gets$ \text{random pop of size $N_0$} 
        \State $\mathbf{y}^g \gets$ \textsc{fobj}($\mathbf{X}^g$) 

        ~    

        ~
        \Repeat
        \State $g \gets g+1$
        \State $\mathbf{X^*} \gets$ \textsc{select}($\mathbf{X}^{g-1}, \mathbf{y}^{g-1}, N$)

        ~

        ~
        \State $\mathbf{X^*} \gets$ \textsc{mate}($\mathbf{X}^*$)
        \State $\mathbf{X}^g \gets$ \textsc{mutate}($\mathbf{X}^*$)
        
        ~

        ~

        ~

        ~

        ~
        \State $\mathbf{y}^g \gets$ \textsc{fobj}($\mathbf{X}^g$) 

        ~

        ~
        \Until stop criteria reached
        \State $\mathbf{X}^{\text{end}} \gets \mathbf{X}^{g}$
        \State \textbf{return} $\mathbf{X}^{\text{end}}$
        \EndProcedure
    \end{algorithmic}
\end{algorithm}
\end{minipage}
\hfill
\begin{minipage}{0.44\textwidth}
\begin{algorithm}[H]
    \scriptsize
    \caption{Proposed ASGA.}\label{algorithm1}
    \hspace*{\algorithmicindent} \textbf{Input}:\\
    \hspace*{\algorithmicindent}\hspace*{10pt} initial population size $N_0$\\
    \hspace*{\algorithmicindent}\hspace*{10pt} population size $N$\\
    \hspace*{\algorithmicindent}\hspace*{10pt} \red{active dimension $M$}\\
    \hspace*{\algorithmicindent}\hspace*{10pt} \red{number backward $B$}\\
    \hspace*{\algorithmicindent}\hspace*{10pt} selection routine \textsc{select}\\
    \hspace*{\algorithmicindent}\hspace*{10pt} mutation routine \textsc{mutate}\\
    \hspace*{\algorithmicindent}\hspace*{10pt} mate routine \textsc{mate}\\
    \hspace*{\algorithmicindent}\hspace*{10pt} objective function \textsc{fobj}\\
    \hspace*{\algorithmicindent}\hspace*{10pt} stop criteria\\
    \hspace*{\algorithmicindent} \textbf{Output}:\\ 
    \hspace*{\algorithmicindent}\hspace*{10pt} final population $\mathbf{X}^{\text{end}}$\\
    \begin{algorithmic}[1]
        \Procedure{ASGA}{}
        \State $g \gets$ 0 
        \State $\mathbf{X}^g \gets$ \text{random pop of size $N_0$} 
        \State $\mathbf{y}^g \gets$ \textsc{fobj}($\mathbf{X}^g$) 
        \State \red{$\mathbf{X}^{\text{AS}} \gets \mathbf{X}^g$}
        \State \red{$\mathbf{y}^{\text{AS}} \gets \mathbf{y}^g$}
        \Repeat
        \State $g \gets g+1$
        \State $\mathbf{X^*} \gets$ \textsc{select}($\mathbf{X}^{g-1}, \mathbf{y}^{g-1}, \red{\textstyle\frac{N}{B}}$)
        \State \red{build AS($\mathbf{X}^{\text{AS}}, \mathbf{y}^{\text{AS}}, M$)}
        \State \red{$\mathbf{X}_M^* \gets$ \textsc{forward}($\mathbf{X}^*$)}
        \State \red{$\mathbf{X}_M^*$} $\gets$ \textsc{mate}(\red{$\mathbf{X}_M^*$})
        \State \red{$\mathbf{X}_M^*$} $\gets$ \textsc{mutate}(\red{$\mathbf{X}_M^*$})
        \For {$\mathbf{x}$ in $\mathbf{X}_M^*$}
        \For {$i \gets 1$ to $B$}
        \State $\mathbf{X}^g \gets$ \textsc{backward}($\mathbf{x}$)
        \EndFor 
        \EndFor 
        \State $\mathbf{y}^g \gets$ \textsc{fobj}($\mathbf{X}^g$) 
        \State \red{$\mathbf{X}^{\text{AS}} \gets \mathbf{X}^{\text{AS}} \cup \mathbf{X}^g$}
        \State \red{$\mathbf{y}^{\text{AS}} \gets \mathbf{y}^{\text{AS}} \cup \mathbf{y}^g$}
        \Until stop criteria reached
        \State $\mathbf{X}^{\text{end}} \gets \mathbf{X}^{g}$
        \State \textbf{return} $\mathbf{X}^{\text{end}}$
        \EndProcedure
    \end{algorithmic}
\end{algorithm}
\end{minipage}\\
 \\

\section{Numerical results}
\label{sec:results}
In this section we are going to present the results obtained by applying the proposed
algorithm, firstly to some test functions that are usually used as benchmarks
for optimization problems. Since this method is particularly suited for
high-dimensional functions, we analyze the optimization convergence for three
different input dimension ($2$, $15$, and $40$), i.e. the number of genes of
each individual.
The second test case we propose is instead a typical engineering problem,
where we optimize the lift-to-drag coefficient of a NACA airfoil which is
deformed using a map $\mathcal{M}: \mathbb{R}^{10} \to \mathbb{R}$ defined in
\cref{sec:naca}. \RA{In this example we opted for the use of a
  surrogate model only to evaluate the individuals' fitness for
computational considerations, since we just want to compare ASGA with
GA. We do not rely on the surrogate for the gradients
approximation. In~\cite{demo2021hull}, instead, we apply
ASGA on a naval engineering hydrodynamics problem, where we do not
rely on a surrogate model of the target function, but instead we
exploit data-driven model order reduction methods to reconstruct the
fields of interest and then compute the function to optimize.}

In both the test cases, in order
to collect a fair comparison, we adopt the same routines for the selection,
the mutation and the crossover steps.
In particular:
\begin{itemize}
\item for the {\it mate} we use the blend BLX-alpha
crossover~\cite{eshelman1993real} with $\alpha= 1.0$, with a mate probability
of $50\%$. With this method, the offspring results:
\begin{equation}
\label{eq:mating_general}
\begin{cases}
\RB{\mathbf{x}_a^i} = (1 - \gamma)\RB{\mathbf{x}_a^{i-1}} + \gamma \RB{\mathbf{x}_b^{i-1}}\\
\RB{\mathbf{x}_b^i} = \gamma \RB{\mathbf{x}_a^{i-1}} + (1 - \gamma)\RB{\mathbf{x}_b^{i-1}}
\end{cases}\quad \text{for}\;\RB{a, b} = 1, \dotsc, N,
\end{equation}
where $\RB{\mathbf{x}_a^{i-1}}$ and $\RB{\mathbf{x}_b^{i-1}}$ refer
       to the parent individuals (\RB{at the $i-1$th generation}), $\RB{\mathbf{x}_a^i}$
       and $\RB{\mathbf{x}_b^i}$ are the mated individuals, \RB{$N$ is the cardinality of population}, and
$\gamma$ is a random variable chosen in the interval $[-\alpha,
1+\alpha)$. \RA{We mention that \cref{eq:mating_general} can recover the graphical
  description of mating in \cref{fig:crossover} if $\gamma$ is taken to
  be discrete, either $0$ or $1$, and applied component-wise.}

\item for the {\it mutation}, a Gaussian
operator~\cite{hinterding1995gaussian} has been used with a mutation
probability of $50\%$. This strategy changes genes by adding a normal noise.
\RA{Since we do not have any knowledge about the low-dimensional space, tuning
the variance of such mutation may result in a not trivial procedure.  This
quantity has in fact to be set in order to explore the input space but,
as the same time, producing minimal differences between parents and
offspring. A fixed variance for both the spaces may cause a too big distance
--- in $l_2$ sense --- between parents and offspring, inhibiting the
convergence.} 
To overcome this potential problem, we correlate the
Gaussian variance with the genes theirselves, ensuring a reasonable mutation in
both spaces. The adopted mutation method is:
\begin{equation}
	\RB{\mathbf{x}_a^i} = \RB{\mathbf{x}_a^{i-1}} + \varepsilon \, \RB{\mathbf{x}}_a^{i-1} \qquad \text{for}\;\RB{a} = 1, \dots, N,
\end{equation}
where $\varepsilon$ is a random variable with probability distribution
$\mathcal{N}(\mu, \sigma^2)$, that is $\varepsilon \sim \mathcal{N}(\mu,
\sigma^2)$, with $\mu = 0$ and $\sigma^2 = 0.1$.
\end{itemize}
Regarding the {\it selection}, since the limited number of individuals per
population, we adopt one of the simplest criteria, by selecting the $N$ best
individuals in terms of fitness.

We also keep fixed the additional
parameters for the AS extension: the number of active dimensions $M$ is set to
$1$, while the number of back-mapped points is $2$. \RB{All the gradients
computations are done using local linear models}. For the actual
computation regarding AS we used the ATHENA\footnote{Freely available
  at \url{https://github.com/mathLab/ATHENA}.} Python 
package~\cite{romor2020athena}. The only varying parameters 
are the size of the initial population $N_0$, the size of population during the
evolution $N$, and the number of generations in the evolutive loop, which are
chosen empirically based on the objective function. We \RA{emphasize} that, due to
the stochastic nature of these methods, we repeated the tests $15$ times, with
different initial configurations, presenting the mean value, the minimum and
the maximum over the $15$ runs.

\subsection{Benchmark test functions}
\begin{figure}
    \centering
    \begin{subfigure}[b]{0.32\textwidth}
	    \includegraphics[width=\textwidth]{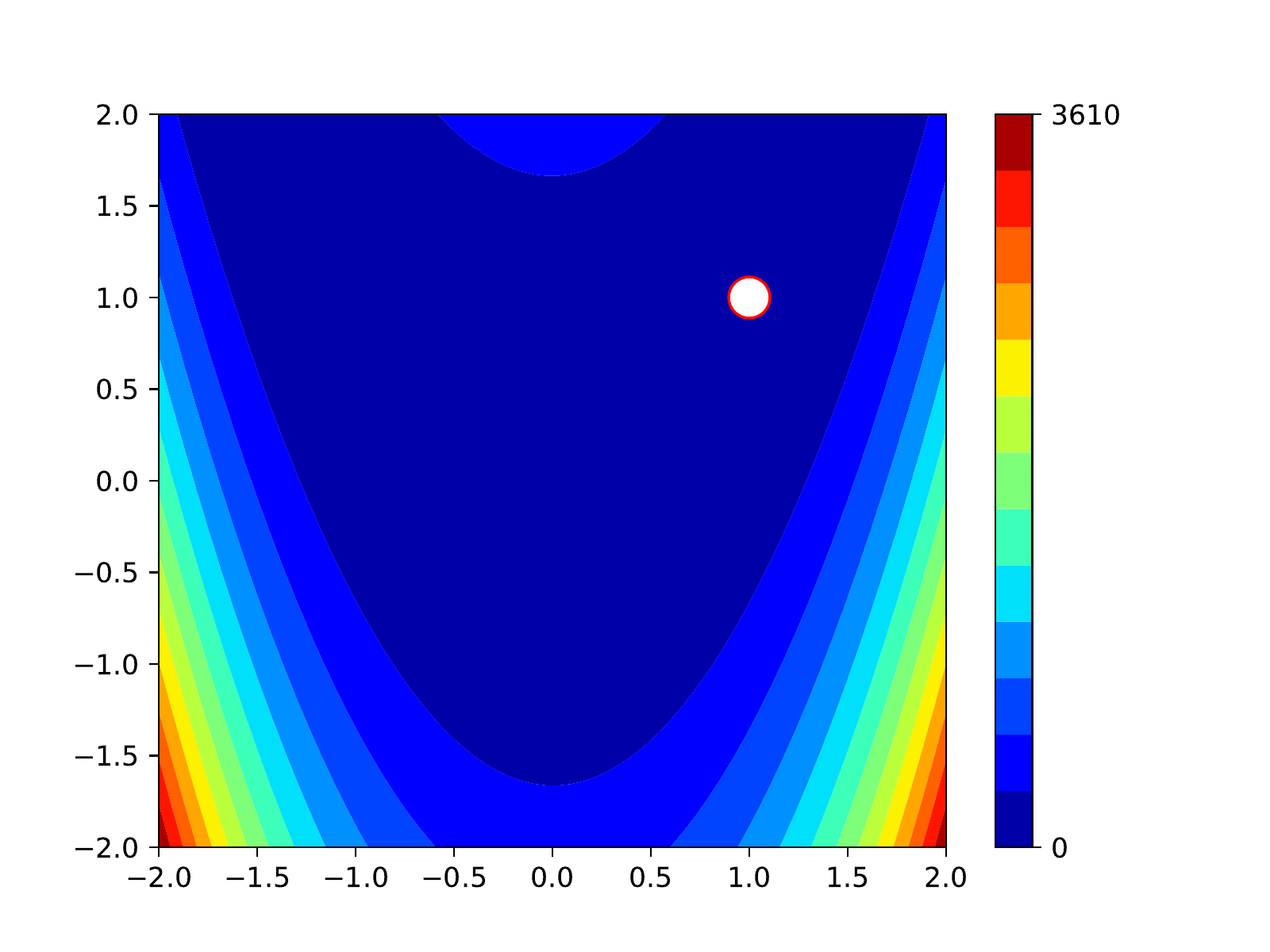}
        \caption{Rosenbrock function.}
        \label{fig:rosenbrock}
    \end{subfigure}
    \begin{subfigure}[b]{0.32\textwidth}
	    \includegraphics[width=\textwidth]{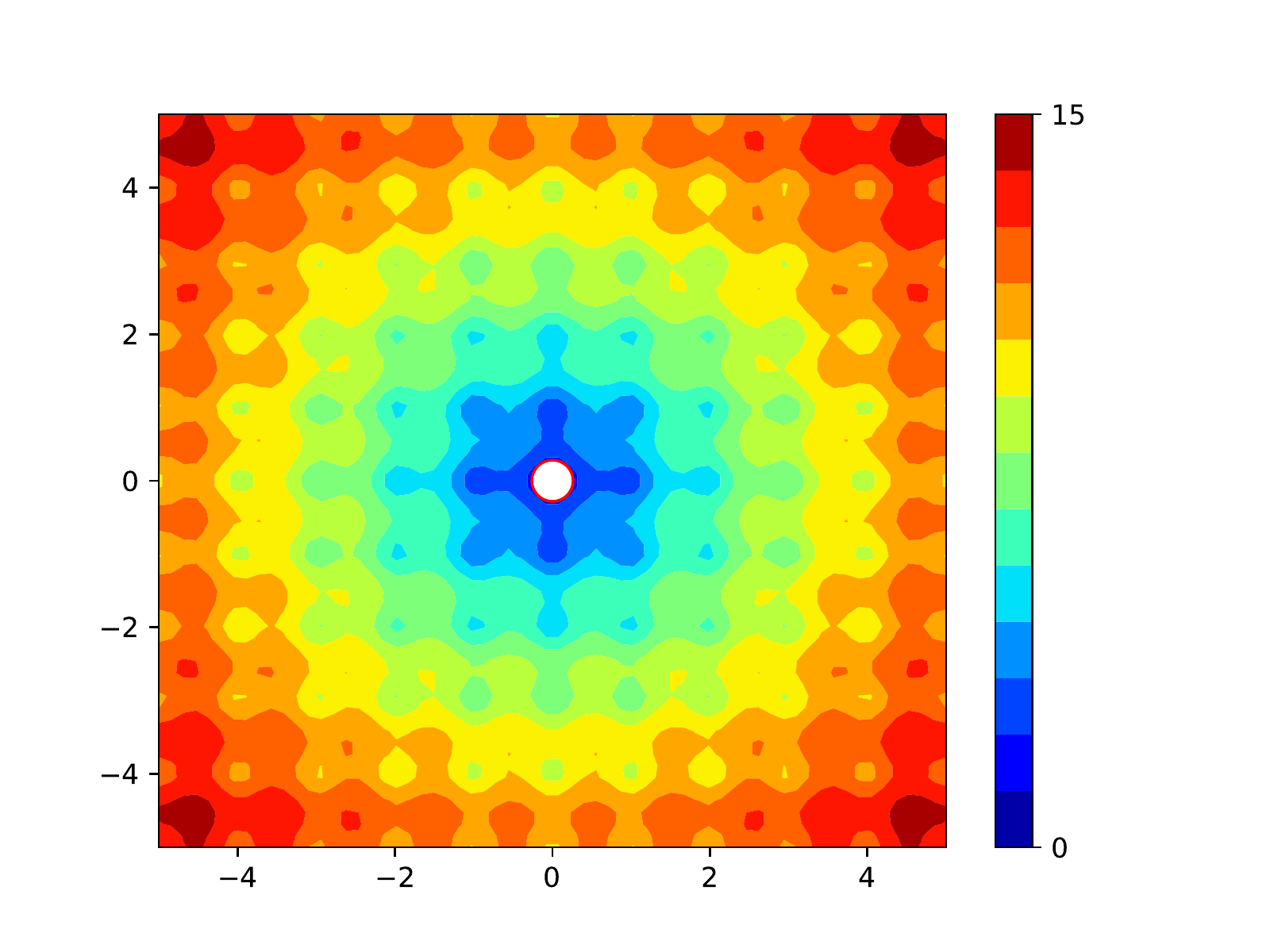}
        \caption{Ackley function.}
        \label{fig:ackley}
    \end{subfigure}
    \begin{subfigure}[b]{0.32\textwidth}
	\includegraphics[width=\textwidth]{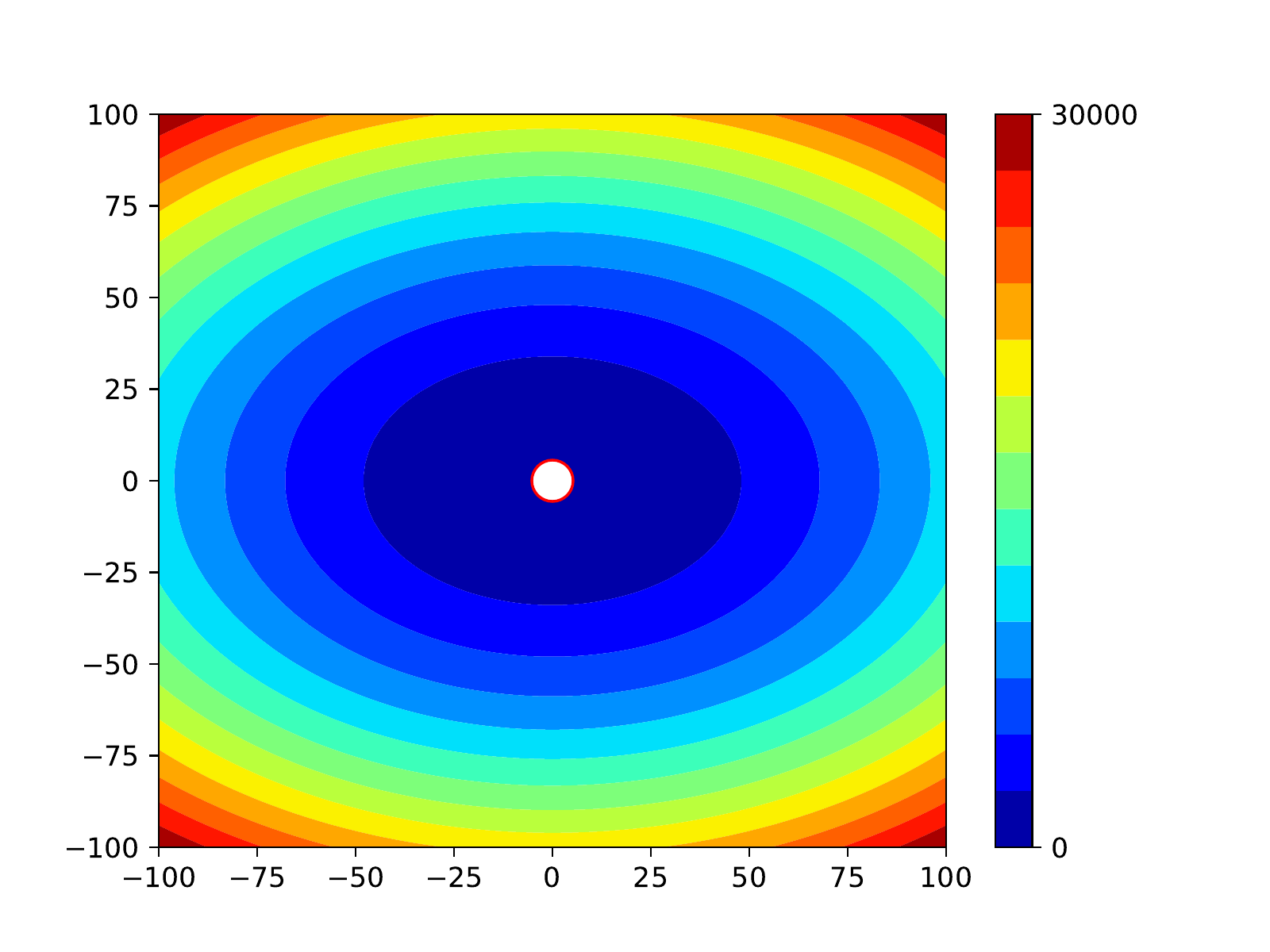}
        \caption{Bohachevsky function.}
        \label{fig:bohachevsky}
    \end{subfigure}

    \begin{subfigure}[b]{0.32\textwidth}
	\includegraphics[width=\textwidth]{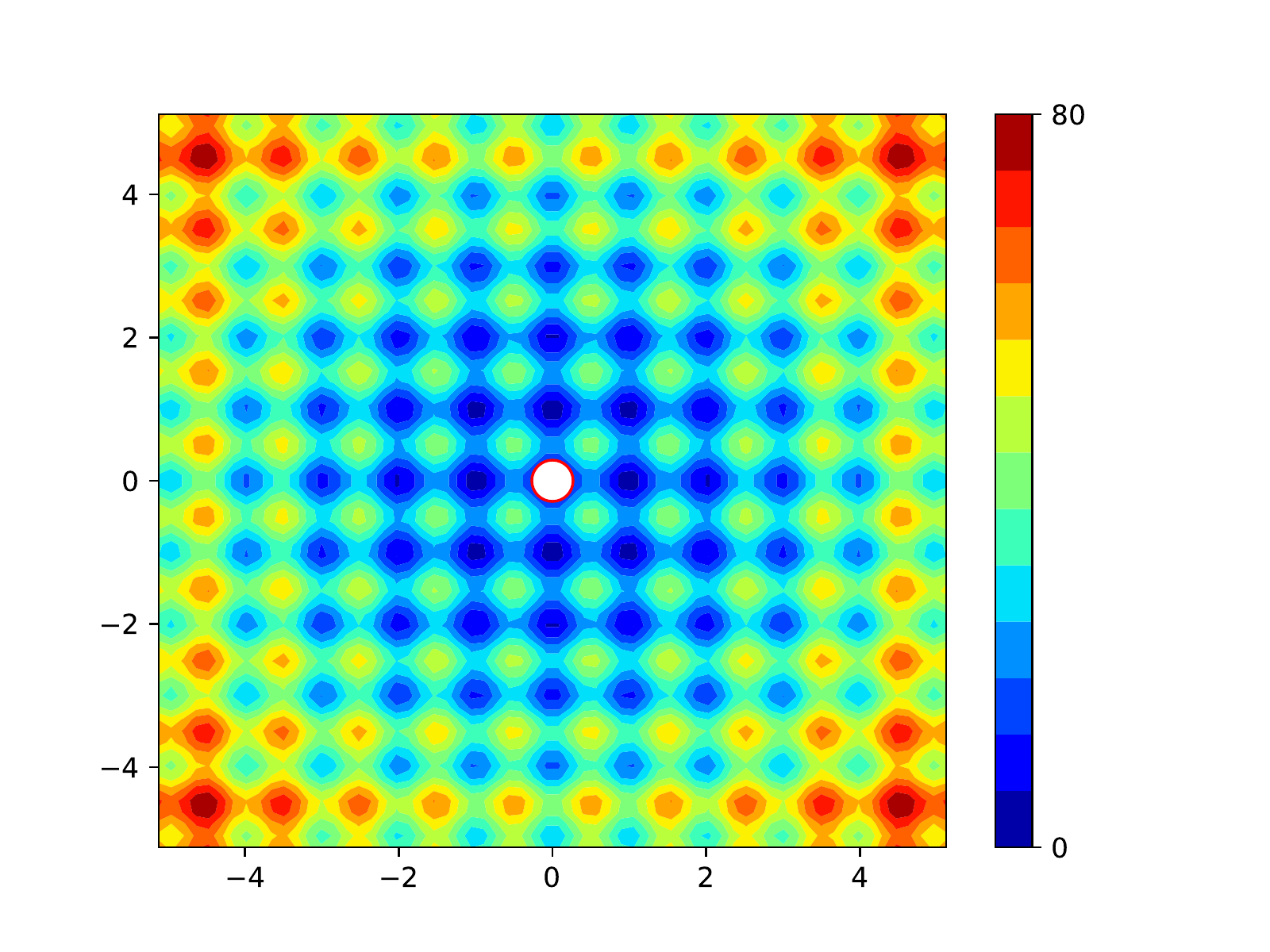}
        \caption{Rastrigin function.}
        \label{fig:rastrigin}
    \end{subfigure}
    \begin{subfigure}[b]{0.32\textwidth}
	    \includegraphics[width=\textwidth]{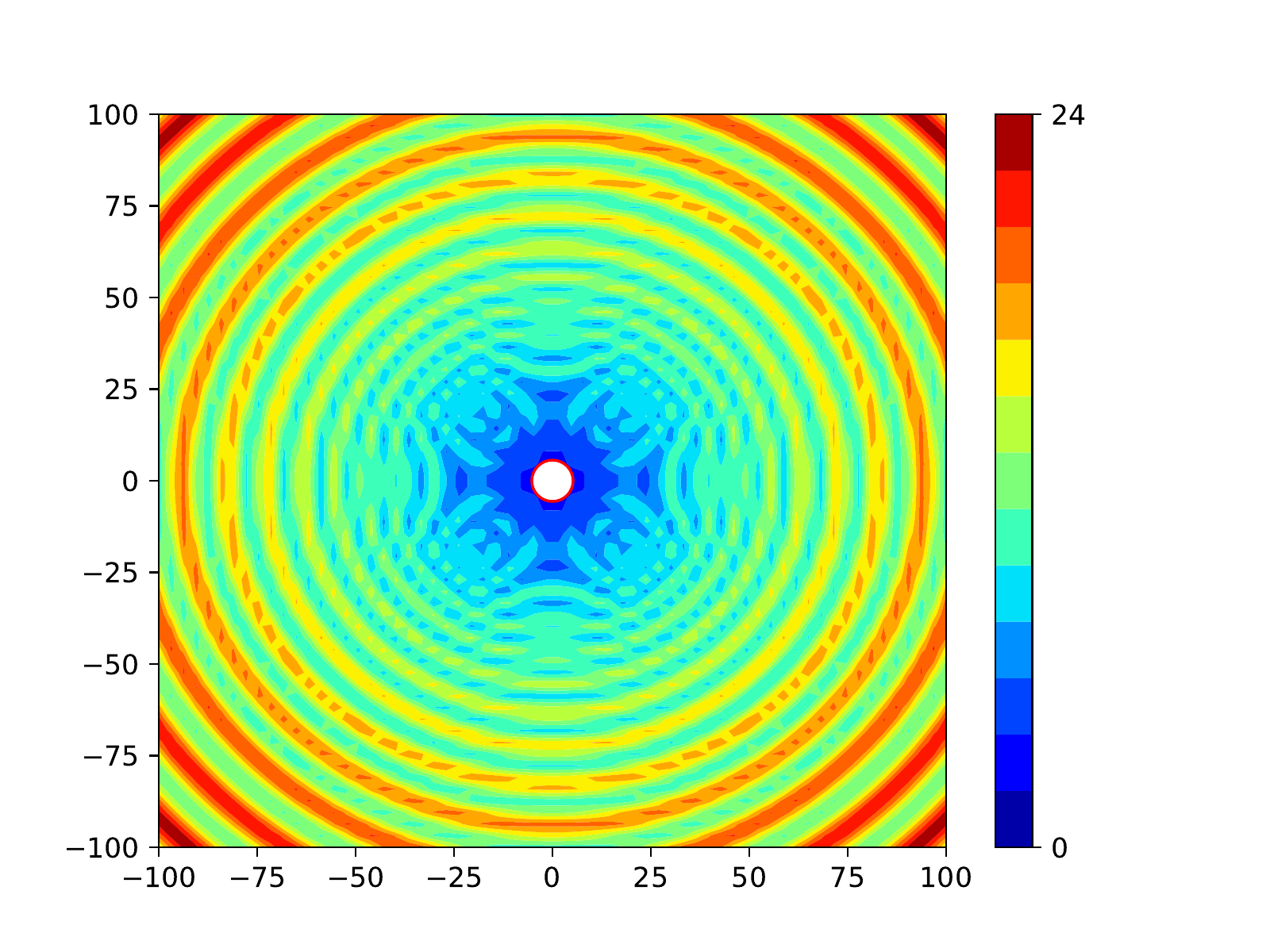}
        \caption{Schaffer N.\ 7 function.}
        \label{fig:schaffer}
    \end{subfigure}
    \begin{subfigure}[b]{0.32\textwidth}
	    \includegraphics[width=\textwidth]{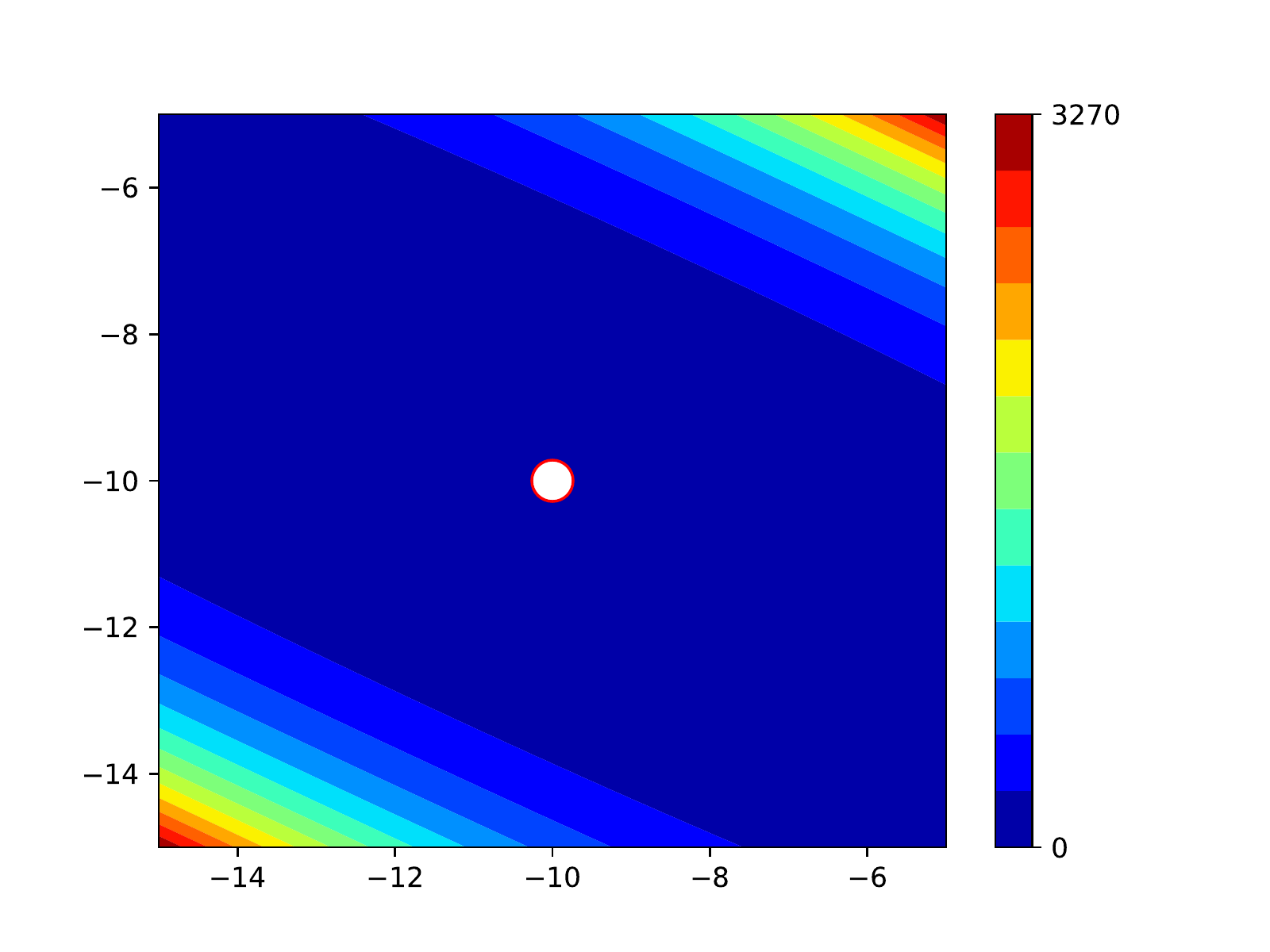}
        \caption{Zakharov function.}
        \label{fig:zakharov}
    \end{subfigure}
    \caption{Benchmark test functions representation in 2D. White dots indicate the global minima.}\label{fig:funcs}
\end{figure}

We applied the optimization algorithm to $6$ different $n$-dimensional test
functions, which have been chosen to cover a large variety of possible shapes.
For all the functions, the results of the proposed method are compared to the
results obtained using the standard genetic approach.  In detail, the
functions we tested are the so called: {\it Rosenbrock}, {\it Ackley}, {\it
Bohachevsky}, {\it Rastrigin}, {\it Schaffer N.\ 7} and {\it Zakharov}. In
Figure~\ref{fig:funcs} we depict the test functions, in their
two-dimensional form. In the
following paragraphs we briefly introduce them before presenting the
obtained results. For a complete literature survey on benchmark functions
for global optimization problems we suggest~\cite{jamil2013survey}.

\paragraph{(a) Rosenbrock function}
The Rosenbrock function is a widespread test function in the context
of global
optimization~\cite{dixon1978global,bect2012sequential,picheny2013benchmark}. 
We choose it as representative of the valley-shaped test functions. The
general $d$-dimensional formulation is the following: 
\begin{equation}
\label{eq:rosenbrock}
f(x) = \sum_{i=1}^{d-1} [100 (x_{i+1} - x_i^2)^2 + (x_i - 1)^2].
\end{equation}
Its global minimum is $f(x^*) = 0$, at $x^* = (1, 1, \dots, 1)$.
As we can see from \cref{fig:rosenbrock} the minimum lies on a
easy to find parabolic valley, but the convergence to the actual
minimum is notoriously difficult. We evaluated the function in the
hypercube $[-5, 10]^d$.

\paragraph{(b) Ackley function}
The Ackley function is characterized by many local minima making it
difficult to find the global minimum, especially for hillclimbing
algorithms~\cite{back1996evolutionary,adorio2005mvf}. The
general $d$-dimensional formulation is the following: 
\begin{equation}
\label{eq:ackley}
f(x) = -a \exp \left( -b  \sqrt{\frac{1}{d}\sum_{i=1}^d x_i^2} \right)
- \exp \left( \sqrt{\frac{1}{d}\sum_{i=1}^d \cos(c x_i)} \right) + a + \exp(1),
\end{equation}
where $a$, $b$, and $c$ are set to $20$, $0.2$, and $2\pi$, respectively.
Its global minimum is $f(x^*) = 0$, at $x^* = (0, 0, \dots, 0)$.
As we can see from \cref{fig:ackley} the function is nearly
flat in the outer region, with many local minima, and the global
minimum lies on a hole around the origin. The function has been evaluated in
the domain $[-15, 30]^d$.

\paragraph{(c) Bohachevsky function}
The Bohachevsky function is a representative of the bowl-shaped
functions. There are many variants and we chose the
general $d$-dimensional formulation as the following: 
\begin{equation}
\label{eq:bohachevsky}
f(x) = \sum_{i=1}^{d-1} (x_i^2 + 2x_{i+1}^2 - 0.3 \cos(3 \pi x_i) -
0.4 \cos(4 \pi x_{i+1}) + 0.7) .
\end{equation}
Its global minimum is $f(x^*) = 0$, at $x^* = (0, 0, \dots, 0)$.
As we can see from \cref{fig:bohachevsky} the function has a
clear bowl shape. This function has been evaluated in the domain $[-100, 100]^d$.

\paragraph{(d) Rastrigin function}
The Rastrigin function is another difficult function to deal with for
global optimization with genetic algorithm due to the large search
space and its many local
minima~\cite{muhlenbein1991parallel}. The general $d$-dimensional
formulation is the following:
\begin{equation}
\label{eq:rastrigin}
f(x) = 10 d + \sum_{i=1}^d [x_i^2 - 10 \cos(2 \pi x_i)].
\end{equation}
Its global minimum is $f(x^*) = 0$, at $x^* = (0, 0, \dots, 0)$.
As we can see from \cref{fig:rastrigin} the function is highly
multimodal with local minima regularly distributed. We evaluated this
function in the input domain $[-5.12, 5.12]^d$.

\paragraph{(e) Schaffer N.\ 7 function}
The Schaffer N.\ 7 function~\cite{schaffer1989study} is a stretched V
sine wave. The general $d$-dimensional
formulation is the following:
\begin{equation}
\label{eq:schaffer}
f(x) = \sum_{i=1}^{d-1} (x_i^2 + x_{i+1}^2)^{0.25} \left [ \sin^2 (50
  (x_i^2 + x_{i+1}^2)^{0.10} ) + 1 \right ] .
\end{equation}
Its global minimum is $f(x^*) = 0$, at $x^* = (0, 0, \dots, 0)$.
As we can see from \cref{fig:schaffer} the function presents
many local minima. The optimization has been performed in the hypercube $[-100,
100]^d$.

\paragraph{(f) Zakharov function}
The Zakharov function is a representative of the plate-shaped
functions. It has no local minima and one global minimum. The general
$d$-dimensional formulation is the following (after a shift):
\begin{equation}
\label{eq:zakharov}
	f(x) = \sum_{i=1}^d (x_i + 10)^2 + \left ( \sum_{i=1}^d \frac{i}{2} (x_i + 10)
	\right )^2 + \left ( \sum_{i=1}^d \frac{i}{2} (x_i+10) \right )^4. 
\end{equation}
We \RA{emphasize} that we used a shifted version with global minimum
$f(x^*) = 0$, at $x^* = (-10, -10, \dots, -10)$. This choice is made
to prove that the proposed method is not biased towards minima around
the origin. 
We can see from \cref{fig:zakharov} the function for $d=2$. We
evaluated the Zakharov function in the domain $[-15, 0]^d$. 

\begin{figure}[htb]
\includegraphics[width=\textwidth]{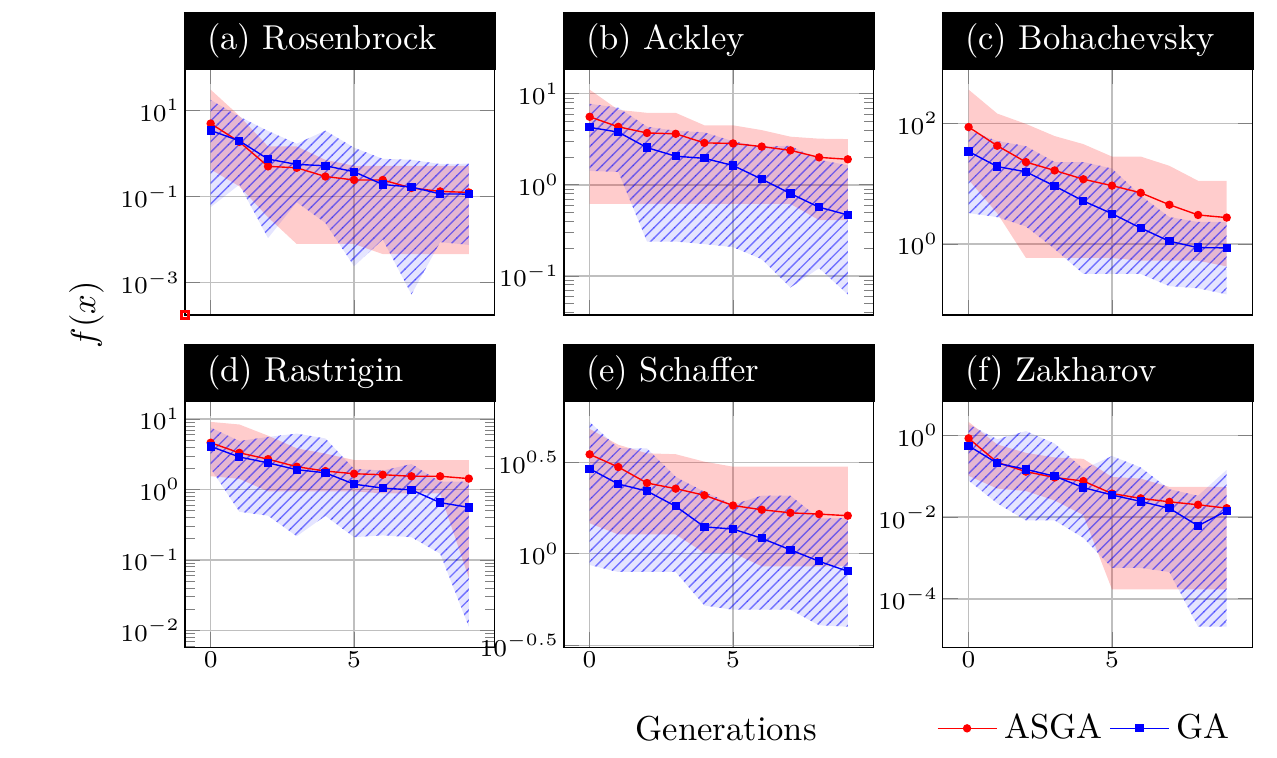}
\caption{Results of the optimization of the benchmark functions in a space of
	dimension $d = 2$. We compare the standard GA (in blue square dots) with
        the proposed algorithm ASGA (in red circle dots) using an initial population of size $30$, while
	the dimension for each generation is fixed to $10$. The solid lines
	represent the mean, over $15$ runs, of the objective function
        corresponding to the best individual at each generation. The shaded areas show
	the interval between minimum and maximum (blue with lines for
        GA, red for ASGA).\label{fig:opt2dim}}
\end{figure}

All the test cases presented share the same hyper-parameters described
at the beginning of this section, except for the population size. For
the $2$-dimensional benchmark functions, the two algorithms are tested
creating $N_0 = 200$ random individuals for the initial population,
then keeping an offspring of dimension $N = 100$. The \cref{fig:opt2dim}
shows the behaviour for all the test functions. For this space dimension, the
two trends are very similar: the usage of the proposed algorithm does not make
the optimization faster, and adds the computational overhead for the AS
construction. Despite that, the results after $10$ generations are
very similar, and we can consider this as a worst case scenario, where
a clear reduction in the parameter space is not possible.

\begin{figure}[htb]
\includegraphics[width=\textwidth]{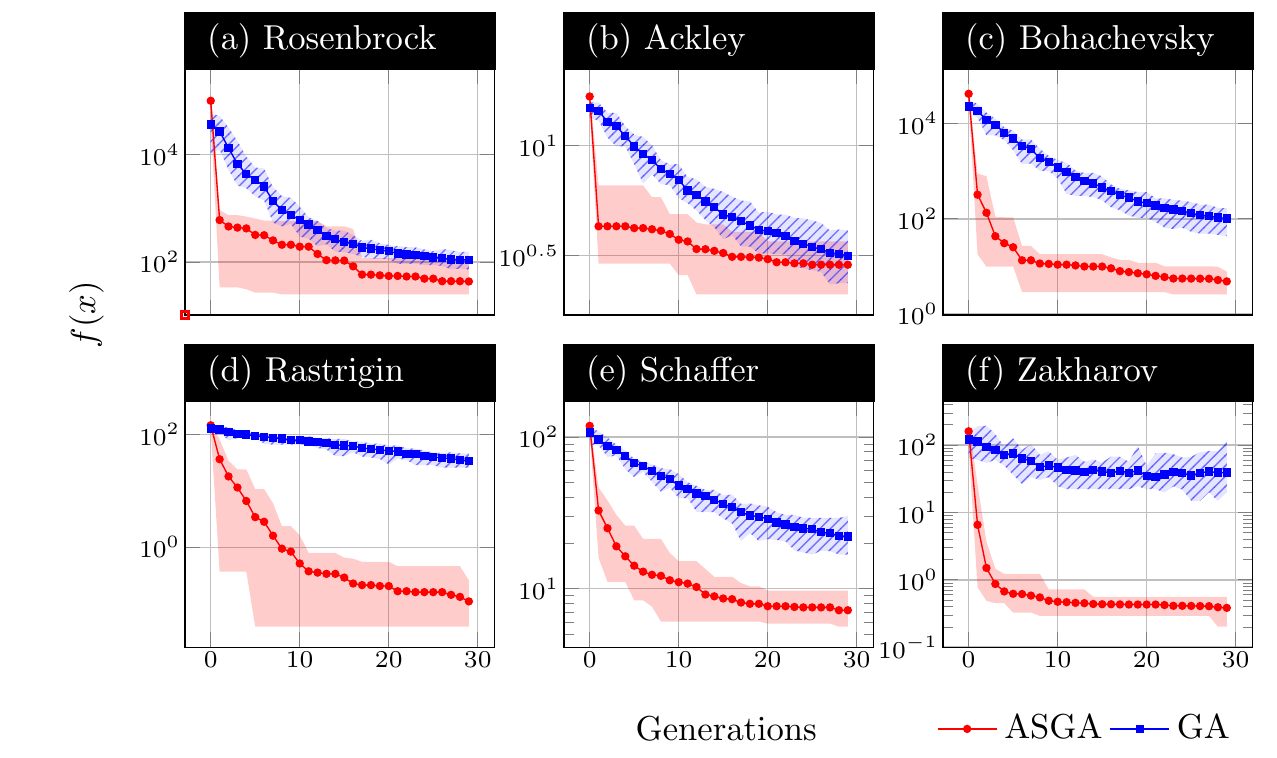}
\caption{Results of the optimization of the benchmark functions in a space of
	dimension $d = 15$. We compare the standard GA (in blue square dots) with
        the proposed algorithm ASGA (in red circle dots) using an initial population of size $2000$, while
	the dimension for each generation is fixed to $200$. The solid lines
	represent the mean, over $15$ runs, of the objective function
        corresponding to the best individual at each generation. The shaded areas show
	the interval between minimum and maximum (blue with lines for
        GA, red for ASGA).\label{fig:opt15dim}}
\end{figure}

The ASGA performance gain changes drastically increasing the number of
dimension to $d = 15$, as demonstrated in \cref{fig:opt15dim}. For such
dimension, the two parameters $N_0$ and $N$ are set to $2000$ and $200$,
respectively.  Starting from this dimension, it is \RA{possible} to note a
remarkable difference between the standard method and the proposed one. The
greater the input dimension, the greater the gain produced by ASGA,
due to the exploitation of the AS reduction. All the benchmarks show a
faster decay, but we can isolate two different patterns in the evolution:
\textit{Rosenbrock} and \textit{Ackley} show a very steep trend in the first
generation gain, while for the next generations the population is not able to
decrease its fitness as much as before, showing a quasi-constant behaviour. The
difference with the standard genetic algorithm is maximized in the first
generation, but even if the evolution using ASGA is not so effective after the
first generation in these two cases, the proposed method is able to achieve
anyway a better result (on average) after $30$ generations. The other
benchmarks instead show a much smoother decay, gradually converging to the optimum.
Despite the lack of the initial step, for these benchmarks the gain with
respect to the standard approach becomes bigger, even if after several
generations the convergence rate decreases.
\begin{figure}[htb]
\includegraphics[width=\textwidth]{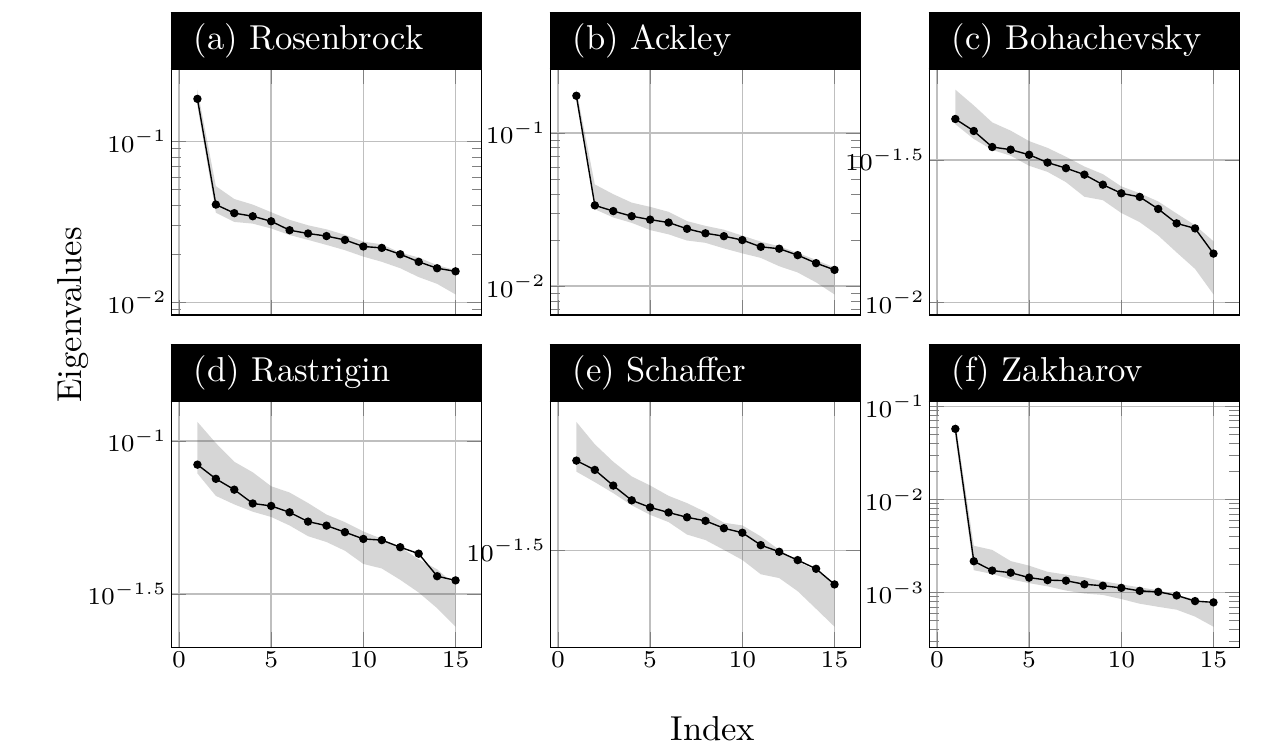}
	\caption{Eigenvalues estimates of the matrix $\mathbf{C}$ in
          \cref{eq:cov_matrix} for all the benchmarks, at
          the first generation, for $d = 15$. The black dots in the plot indicate
          the eigenvalues, while the grey area is defined by the
          bootstrap intervals.\label{fig:as_eigs}}
\end{figure}
In order to better understand these differences, we investigate the spectra of the AS
covariance matrices for all the benchmarks, which are reported in
\cref{fig:as_eigs}. The patterns individuated in the optimizations are
partially reflected in the eigenvalues: \textit{Rosenbrock}, \textit{Ackley}
and \textit{Zakharov} have an evident gap between the first and the second
eigenvalues, which results in a better approximation (of the original function) in
the $1$-dimensional subspace. However, the order of magnitude of the first
eigenvalue is different between the three functions: for \textit{Rosenbrock}
and \textit{Ackley} the magnitude is greater than~$\num{0.1}$ whereas for
\textit{Zakharov} is around~$\num{0.05}$.

\begin{figure}[p!]
\includegraphics[width=\textwidth]{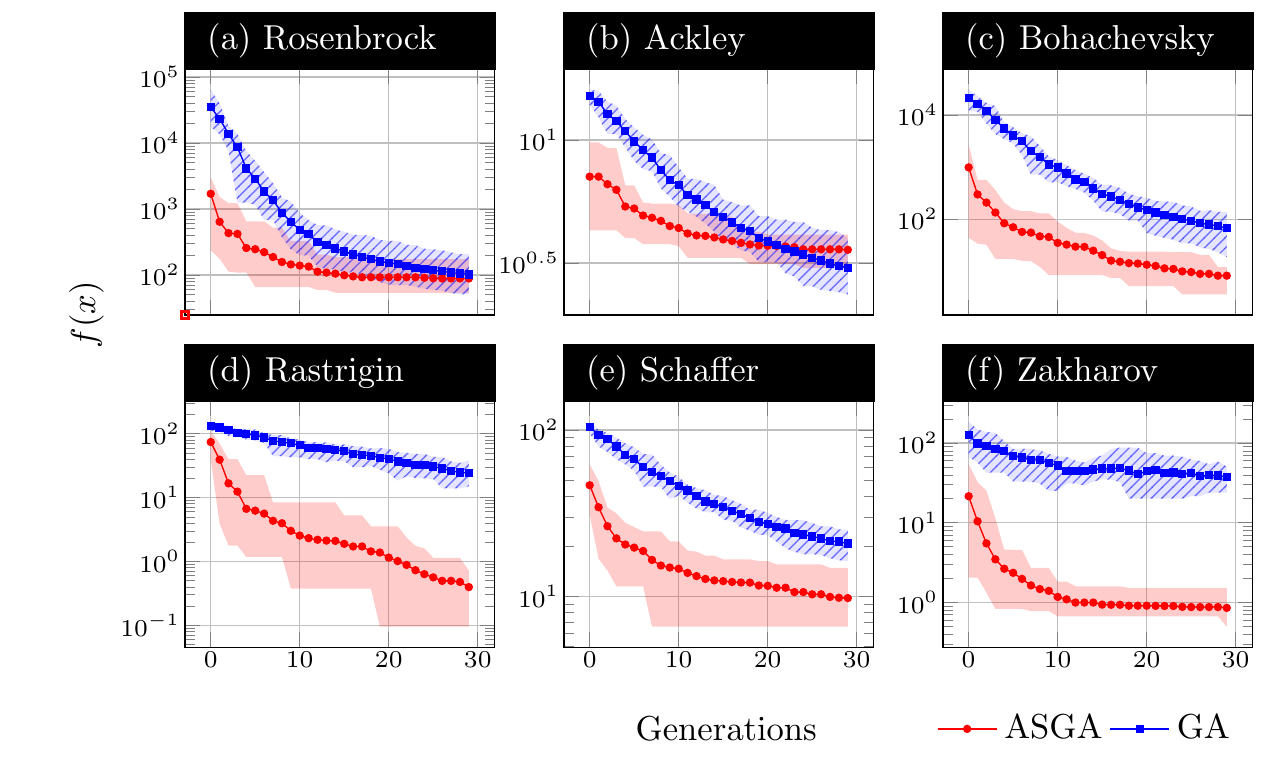}
\includegraphics[width=\textwidth]{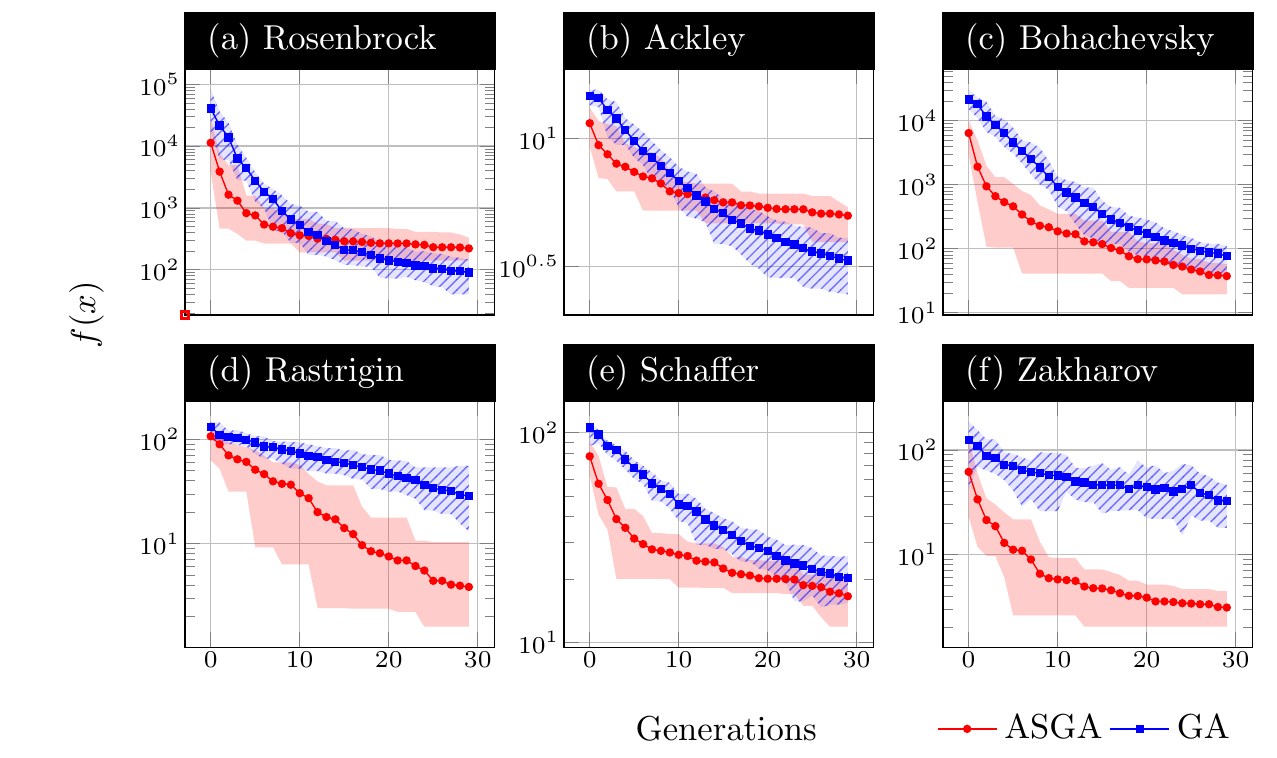}
\caption{\RA{Results of the optimization of the benchmark functions in a space of
	dimension $d = 15$, with active dimension $M=2$ (top) and $M=5$
	(bottom). We compare the standard GA (in blue square dots) with the
	proposed algorithm ASGA (in red circle dots) using an initial
	population of size $2000$, while the dimension for each generation is
	fixed to $200$. The solid lines represent the mean, over $15$ runs, of
	the objective function
        corresponding to the best individual at each generation. The shaded areas show
	the interval between minimum and maximum (blue with lines for
        GA, red for ASGA).}\label{fig:opt15dim_as}}
\end{figure}

\RA{Since for all the tests, the ASGA approach performs better than its
classical counterpart despite the absence, in some cases, of an
evident spectral gap in the AS covariance matrix, we perform further
investigations. In particular, we use the
same tests as before ($15$-dimensional functions) but with a different
number of active dimensions, i.e.\ $M =
\{2, 5\}$, instead of $M = 1$. In \cref{fig:opt15dim_as} we show the
comparison between the classical GA and the ASGA outcomes. It is
possible to note that by increasing the active
dimension, the differences between the performances of the two methods become
smaller. Only for the Rosenbrock and for the Ackley functions we can see that ASGA with $M = 5$
is not able to reach the same order of magnitude reached by GA (we remark the original
space has dimension $d=15$).}

\begin{figure}[htb]
\includegraphics[width=\textwidth]{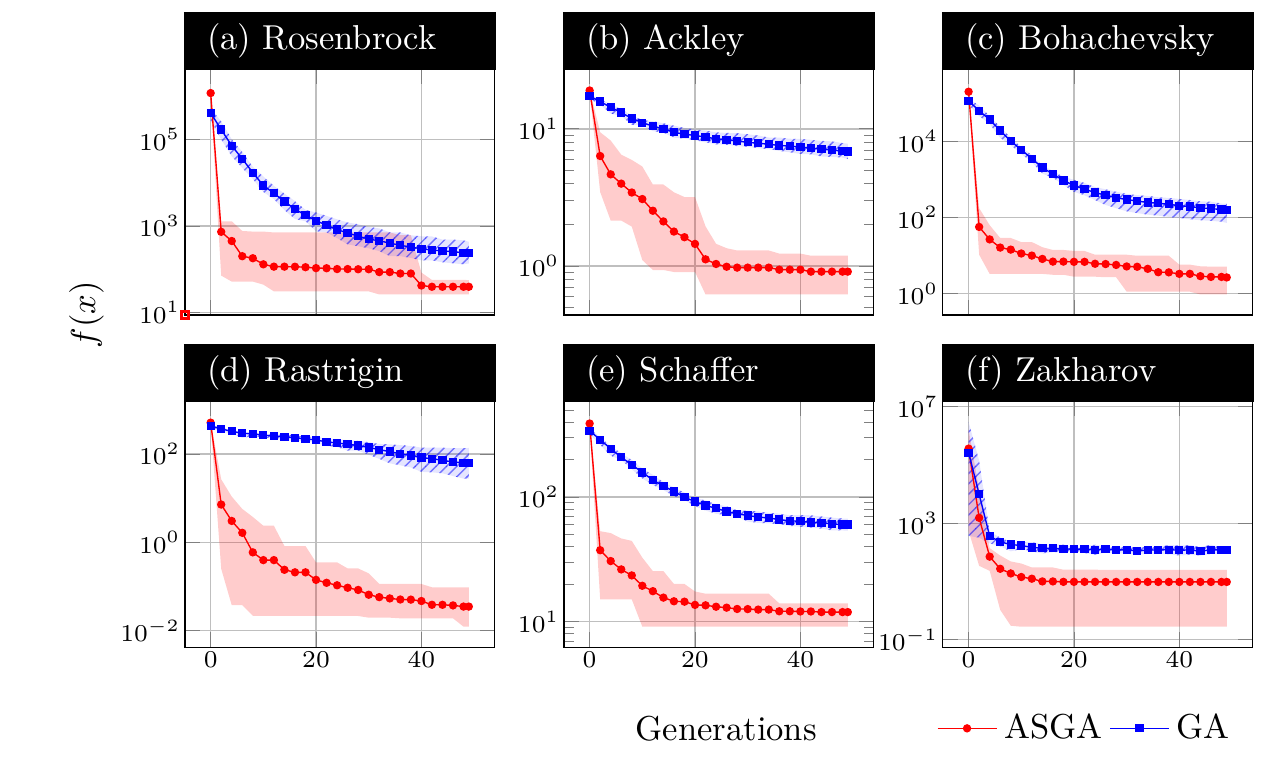}
\caption{Results of the optimization of the benchmark functions in a space of
	dimension $d = 40$. We compare the standard GA (in blue square dots) with
        the proposed algorithm ASGA (in red circle dots) using an initial population of size $5000$, while
	the dimension for each generation is fixed to $1000$. The solid lines
	represent the mean, over $15$ runs, of the objective function
        corresponding to the best individual at each generation. The shaded areas show
	the interval between minimum and maximum (blue with lines for
        GA, red for ASGA).\label{fig:opt40dim}}
\end{figure}

Increasing the input dimension to $d = 40$ \RA{shows a much clearer benefit} in using
the proposed method, as we can see in \cref{fig:opt40dim}. Here we set $N_0
= 5000$ and $N = 1000$. \RA{We specify that we set the active dimension $M=1$}.
Also with this dimensionality, we are able to isolate
two main behaviours in the convergence of the six benchmarks: a very steep
trend in the first generation, and a more smooth one, but still equally
effective. The interesting thing is that some benchmarks do not reflect the
behaviour collected with $d = 15$. While \textit{Rosenbrock},
\textit{Rastrigin} and \textit{Zakharov} show a similar convergence
rate for ASGA, the other benchmarks present a change in the slope. The
different behaviours observed for the same benchmarks 
evaluated at different input dimensional spaces is due to the fact
that the method is sensitive to the approximation accuracy of the
gradients of the model function with respect to the input data. This
is an issue inherited by the application of AS. Moreover, since we are
keeping just one active variable, we are discarding several
information, thus the representation of the function along the active
subspace could present some noise. So the genetic procedure
enhanced by AS is able to converge fast to the optimum, but this
optimum may be --- for the space simplification --- distant to the true optimum.
\RA{From the tests with higher active dimension, we note that the improvement
in the first iterations is not as rapid as by using a $1$-dimensional AS. Also, keeping
more active dimensions, the performance of ASGA becomes similar to the
standard GA. We can conclude that with one (or few) active dimension,
ASGA reaches the global minimum with less function evaluations, but we
get stuck with the projection error introduced by AS, whereas by
increasing the active dimension we reduce the projection error but we
lose the effectiveness of the evolution steps in a reduced space.} A possible solution for
this problem can be a smarter (and dynamical) strategy to select the number of
active dimensions.

Over all the three test cases, where we vary the input space dimension, the
performance of ASGA are better or equal than the standard GA. In
\cref{tab:bench} we summarize the relative gain on average achieved after the entire evolution
and after only one generation, divided by test function, both with the GA and
ASGA methods. The relative gain is computed as the mean over $15$ runs
of the ratio between the objective function
evaluated at the best-fit individual at the beginning of the evolution
($f(x_{\text{opt}}^{0})$) and the
objective function evaluated at the best-fit individual after $k$
generations, with $k=1$ and $k=N_{\text{gen}}$, where $N_{\text{gen}}$
is the maximum number of generations depending on the input
dimension. This relative gain $G(k)$ reads,
\begin{equation}
  \label{eq:gain}
G(k) = \frac{1}{15}\sum_{i=1}^{15} f(x_{\text{opt}_i}^{0}) /  f(x_{\text{opt}_i}^k),
\end{equation}
where $x_{\text{opt}_i}^k$ is the best-fit individual of the population at the $i$-th
run and at $k$-th generation. Highest values correspond to a more effective
optimization, and for dimensions $15$ and $40$ we can see from
\cref{tab:bench} that ASGA performs better than standard GA for all
the benchmarks. Even the gain after just one
evolutive iteration is bigger in all the collected tests, reaching in some
cases some order of magnitude of difference with respect to GA.
These results suggest that despite the computational overhead for the
construction of AS and the back-mapping, an application of the ASGA
over the standard GA produces usually better or at least comparable results for
a fixed generation.

\begin{table}[htb]
\begin{center}
    \caption{Summary comparison between GA and ASGA with respect to
      the gain $G(k)$ defined in \cref{eq:gain}. We compare the
      gain for the first and for the last generation.\label{tab:bench}}
\footnotesize
\begin{tabular}{ ll|rr|rr|rr }
\hline
\hline
  \multirow{2}{*}{\bf Function}  & \multirow{2}{*}{\bf Method} &
      \multicolumn{2}{c}{\bf dim = 2}  &
      \multicolumn{2}{c}{\bf dim = 15} &
      \multicolumn{2}{c}{\bf dim = 40} \\
    && $G(9)$ & $G(1)$ &
      $G(29)$ & $G(1)$ &
      $G(49)$ & $G(1)$ \\
\hline
\rowcolor{Gray}
                         & GA   & 9.17 & 1.13 & 4.71 & 1.03 & 2.53  & 1.10\\
\rowcolor{Gray}
\multirow{-2}{*}{Ackley} & ASGA & 2.93 & 1.29 & 5.81 & 3.89 & 20.91 & 3.00\\

                              & GA   & 39.58 & 1.78 & 223.86  &   1.22 &   729.05 & 1.81   \\
\multirow{-2}{*}{Bohachevsky} & ASGA & 31.66 & 2.04 & 8608.41 & 130.72 & 75104.33 & 3548.70\\

\rowcolor{Gray}
                            & GA   & 7.34 & 1.41 &    3.80 & 1.05 &     6.97 &  1.17\\
\rowcolor{Gray}
\multirow{-2}{*}{Rastrigin} & ASGA & 3.24 & 1.39 & 1343.41 & 4.00 & 14738.40 & 71.77\\

                             &  GA  & 30.04 & 1.74 &  335.33 &   1.34 &  1723.89 &    2.42\\
\multirow{-2}{*}{Rosenbrock} & ASGA & 39.68 & 2.66 & 2343.57 & 167.48 & 29747.56 & 1600.24\\

\rowcolor{Gray}
                           &  GA  & 3.64 & 1.21 &  4.83 & 1.11 &  5.66 &  1.18\\
\rowcolor{Gray}
\multirow{-2}{*}{Schaffer} & ASGA & 2.16 & 1.17 & 16.41 & 3.61 & 32.57 & 10.38\\

                           &  GA  & 38.59 & 2.65 &   3.11 &  1.07 &  2148.39 & 26.50\\
\multirow{-2}{*}{Zakharov} & ASGA & 51.14 & 3.88 & 417.86 & 24.46 & 37739.61 & 237.48 \\
  \hline
  \hline
\end{tabular}
\end{center}
\end{table}

\subsection{Shape design optimization of a NACA airfoil}
\label{sec:naca}
Here we present the shape design optimization of a NACA 4412
airfoil~\cite{abbott2012theory}. Since the purpose of this work is limited to
the extension of the GA, we briefly present the details of the
complete model, with a quick overview of the application. To reproduce
the full order simulations please refer to~\cite{tezzele2019mortech}.

\RB{Let be given the unsteady incompressible Navier-Stokes equations
described in an Eulerian framework on a parametrized space-time domain
$Q(\mupar) = \Omega \times [0,T] \subset
\R^d\times\R^+$, $d=2,3$ with the 
velocity field denoted by $\mathbf{u}: Q(\mupar) \to \R^d$, and the 
pressure field by $p: Q(\mupar) \to \R$, such that:
\begin{equation}
\label{eq:navstokes}
\begin{cases}
\mathbf{u_t}+ \mathbf{\nabla} \cdot (\mathbf{u} \otimes \mathbf{u})- \mathbf{\nabla} \cdot
2 \nu \mathbf{\nabla^s} \mathbf{u}=-\mathbf{\nabla}p &\mbox{ in } Q(\mupar),\\
\mathbf{\nabla} \cdot \mathbf{u}=\mathbf{0} &\mbox{ in } Q(\mupar),\\
\mathbf{u} (t,x) = \mathbf{f}(\mathbf{x}) &\mbox{ on } \Gamma_{\text{in}} \times [0,T],\\
\mathbf{u} (t,x) = \mathbf{0} &\mbox{ on } \Gamma_{0}(\mupar) \times [0,T],\\ 
(\nu\nabla \mathbf{u} - p\mathbf{I})\mathbf{n} = \mathbf{0} &\mbox{ on } \Gamma_{\text{out}} \times [0,T],\\ 
\mathbf{u}(0,\mathbf{x})=\mathbf{k}(\mathbf{x}) &\mbox{ in } Q(\mupar)_0,\\            
\end{cases}
\end{equation}
holds.  Here, we have that $\Gamma = \Gamma_{\text{in}} \cup \Gamma_{0} \cup \Gamma_{\text{out}}$
is the boundary of $\Omega$ and it is composed by the inlet boundary
$\Gamma_{\text{in}}$, the outlet boundary $\Gamma_{\text{out}}$ and the
physical walls $\Gamma_0(\mupar)$. The
term $\mathbf{f}(\mathbf{x})$ stands for the stationary non-homogeneous boundary
condition, whereas $\mathbf{k}(\mathbf{x})$ indicates the initial condition for
the velocity at $t=0$. Shape changes are applied to the boundary
$\Gamma_0(\mupar)$ corresponding to the airfoil 
wall, which in the undeformed configuration corresponds to
the 4-digits, NACA~4412 wing profile. Such shape modifications are
associated to numerical parameters contained in the 
vector $\mupar \in \R^k$ with $k=10$.}

\begin{figure}
\centering
\includegraphics[width=.9\textwidth]{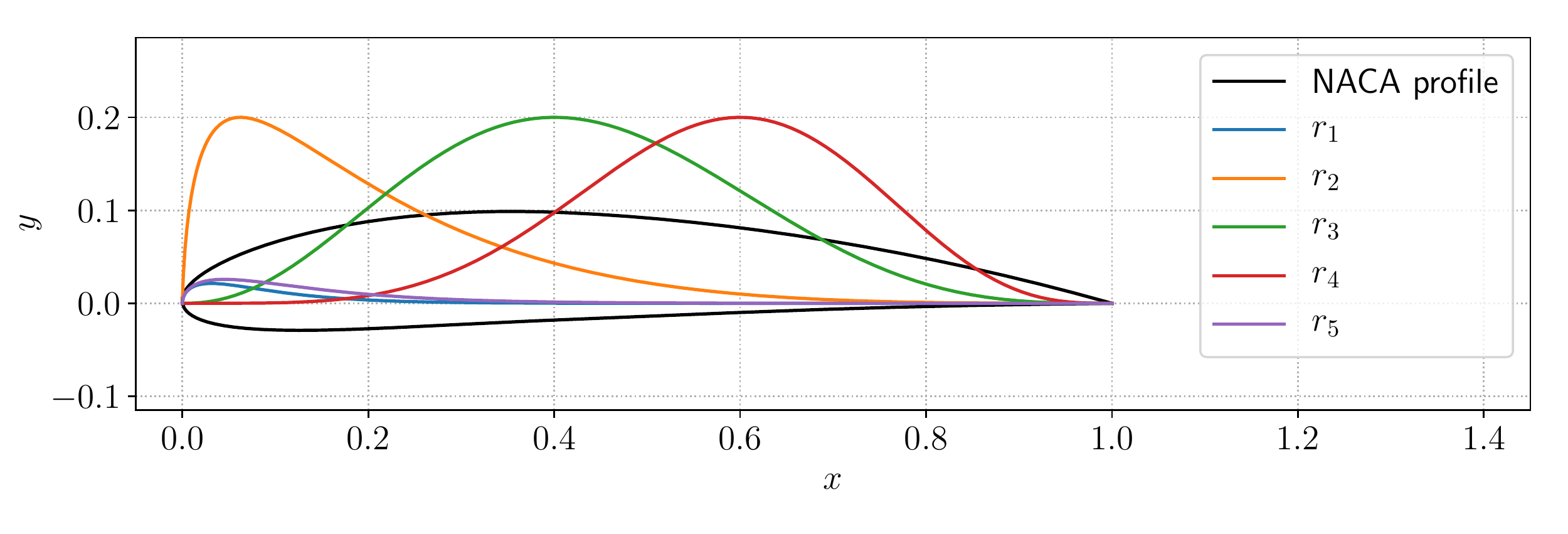}
\caption{NACA 4412 profile with the $5$ shape functions $r_i$ rescaled
  by a factor equal to $0.2$.\label{fig:bumps}}
\end{figure}

As geometrical deformation map
$\mathcal{M}$ we adopt the shape morphing proposed
in~\cite{hicks1978wing}, where 5 shape functions $r_i$ are added to the
upper and lower part of the airfoil profile denoted by $y^+$ and
$y^-$, respectively. Each shape function is multiplied by a possible
different coefficient as in the following
\begin{equation}
  y^+ = \overline{y^+} + \sum_{i=1}^5 a_i r_i , \qquad
  y^- = \overline{y^-} - \sum_{i=1}^5 b_i r_i ,
\end{equation}
where the bar denotes the reference undeformed profile.
These $10$ coefficients ($a_i$ and $b_i$) represent the input
parameters $\mupar \in \mathbb{D} := [0, 0.03]^{10}$. In
\cref{fig:bumps} we depict the NACA 4412 together with the $5$ rescaled
 shape functions $r_i$.
The output function we want to maximize is the lift-to-drag coefficient, one of the
typical quantities of interest in aeronautical problems. To recast the
problem in a minimization setting, we just minimize the opposite of
the coefficient. To compute it, we
model a turbulent flow pasting around the 2D airfoil using the incompressible
Reynolds Averaged Navier--Stokes equations. Regarding the main numerical
settings, we adopt a finite volume approach with the Spalart--Allmaras model,
with a computational grid of $46500$ degrees of freedom. The flow velocity, at
the inlet boundary, is set to $1$ m/s, while the Reynolds number is fixed
to $50000$.
For the detailed problem formulation we refer to the experiments conducted
in~\cite{tezzele2019mortech}.

Instead of running the high-fidelity solver for any new untested parameter, we optimize a
\RB{RBF} response surface built using the initial dataset. Due to the stochastic nature of the method, also in
this test case we test the methods for several initial settings --- $25$
different runs --- making the total computational load very high. Thus, we
decided to build a response surface using a dataset of $333$ samples, computed
with the numerical scheme described above, mimicking at the same time a typical
industrial workflow.

The objective function $f_{\text obj} (\mupar ) : \mathbb{D} \subset
\R^{10} \to \R$ we are going to minimize is the following:
\begin{equation}
f_{\text obj}(\mupar) = \begin{cases}
g(\mupar)\quad\text{if }\, \mupar  \in \mathbb{D},\\
\alpha\quad\quad\,\,\text{if }\, \mupar \notin \mathbb{D},\\
\end{cases}
\end{equation}
where $g(\mupar)$ is the response surface built using the radial basis
function (RBF) interpolation technique~\cite{buhmann2003radial} over the
samples, while $\alpha \in \R$ is a penalty constant. \RA{To prevent
  the evolution from creating} new individuals that do not belong to
$\mathbb{D}$, we impose a penalty factor $\alpha = 10$. We recall that
we minimize the opposite of the lift-to-drag coefficient.

\begin{figure}
\centering
\input{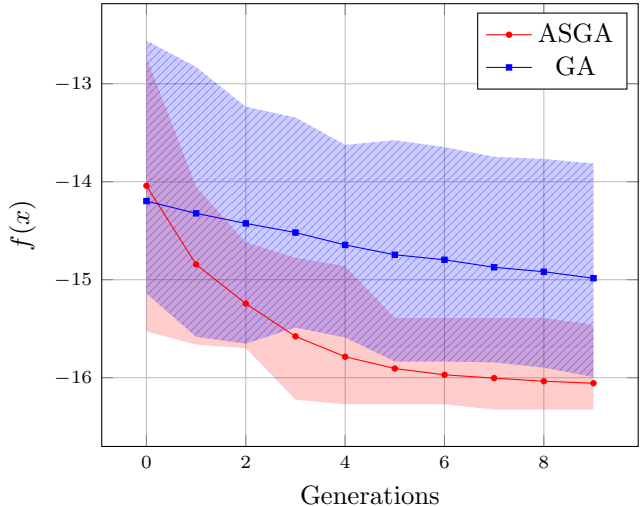}
\caption{Results of the optimization of the NACA airfoil design in a
  $10$-dimensional space. We compare the standard GA (in blue square dots) with
        the proposed algorithm ASGA (in red circle dots) using an initial population of size $20$, while
	the dimension for each generation is fixed to $10$. The solid lines
	represent the mean, over $25$ runs, of the objective function
        corresponding to the best individual at each generation. The shaded areas show
	the interval between minimum and maximum (blue with lines for
        GA, red for ASGA).\label{fig:optnaca}}
\end{figure}

\cref{fig:optnaca} reports the evolution of the best-fit individual
over $10$ generations. Also in this case, we apply the 
proposed algorithm and the standard GA to $25$ different initial settings, using
an initial population size $N_0 = 20$ and selecting at each generations the $N =
10$ best-fit individuals for the offspring.
The plot depicts the mean best-fit individual with solid lines, whereas
the shaded areas show the interval between the minimum and maximum (of the $25$
runs) for each generation. Even if the dimension of the parameter space
is not very high ($10$), we can see that on average the proposed algorithm is
able to converge faster. The difference between the two
methods \RA{is not as remarkable} as in a higher dimensional test case, but we can see
that the best run using standard GA is slightly worse than the mean optimum
achieved by ASGA. \RA{This again demonstrates the value in the
  proposed method.} Moreover, we \RA{emphasize} that also in this case the
decay of the objective function in the first generations with ASGA is faster.

\section{Conclusions}
\label{sec:conclusions}
In this work, we have presented a novel approach for optimization problems
coupling the \RA{supervised learning technique called active subspaces (AS)} with the standard genetic
algorithm (GA). We have \RA{demonstrated} the benefits of such method by applying it
to some benchmark functions and to a realistic engineering
problem. The proposed method achieves faster convergence to the
optimum, since the individuals evolve only along few 
principal directions (discovered exploiting the AS property). Further, from the results it
emerges that the gain induced from the ASGA method is greater for
high-dimensionality functions, making it particularly suited for models with
many input parameters.

This new method can also be integrated in numerical pipelines
involving model order reduction and reduction in parameter
space. Reducing the number of input parameters is a key ingredient to
improve the computational performance and to allow the study of very
complex systems.

Since the number of active dimensions is important for the accuracy of the AS,
future developments will focus on an efficient criterion to select
\RA{dynamically} the number
of AS dimensions, which in the presented results are kept fixed.
Future studies will also address the problem of incorporating non-linear
extensions of active subspaces into the ASGA,
focusing on the construction of a proper back-mapping from the reduced
space to the original full parameter space. 

\appendix
\section{On ASGA convergence}
\label{appendix}
The aim of this section is to provide further insights about the
convergence of the ASGA method. We perform a single run on all the benchmark
functions presented above, in a space of dimension $d = 2$. We kept
unaltered all the ASGA numerical settings described in \cref{sec:results}, so for
all the details we refer to such section. We emphasize that we used the
same hyper-parameters of the $2$-dimensional optimization test, except for the
number of generations which we increased to $100$.

\begin{figure}[htb]
\centering
\begin{minipage}[b]{0.29\textwidth}
\includegraphics[width=\textwidth,trim={0 .5cm 0 1.3cm},clip]{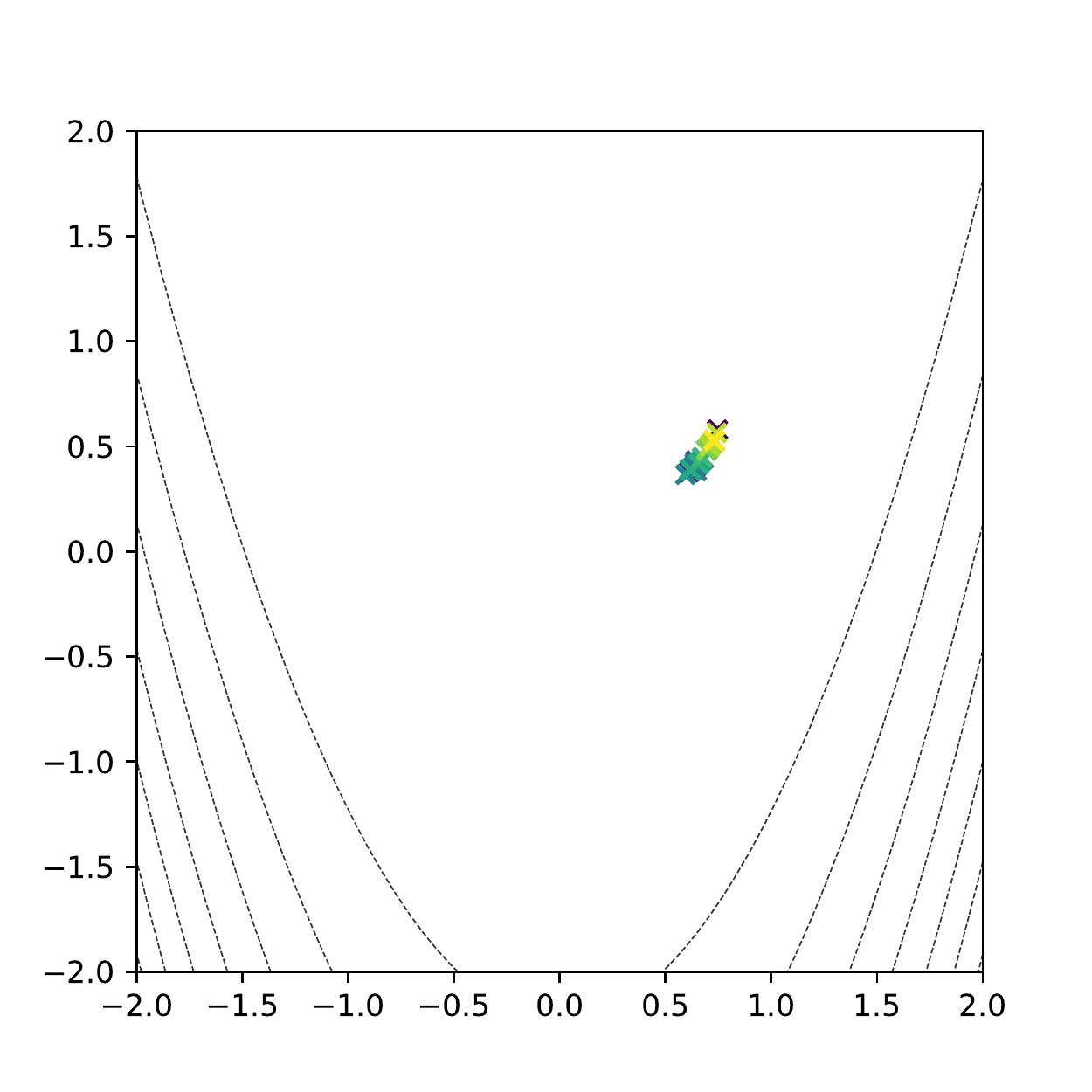}
\includegraphics[width=\textwidth,trim={0 .5cm 0 1.3cm},clip]{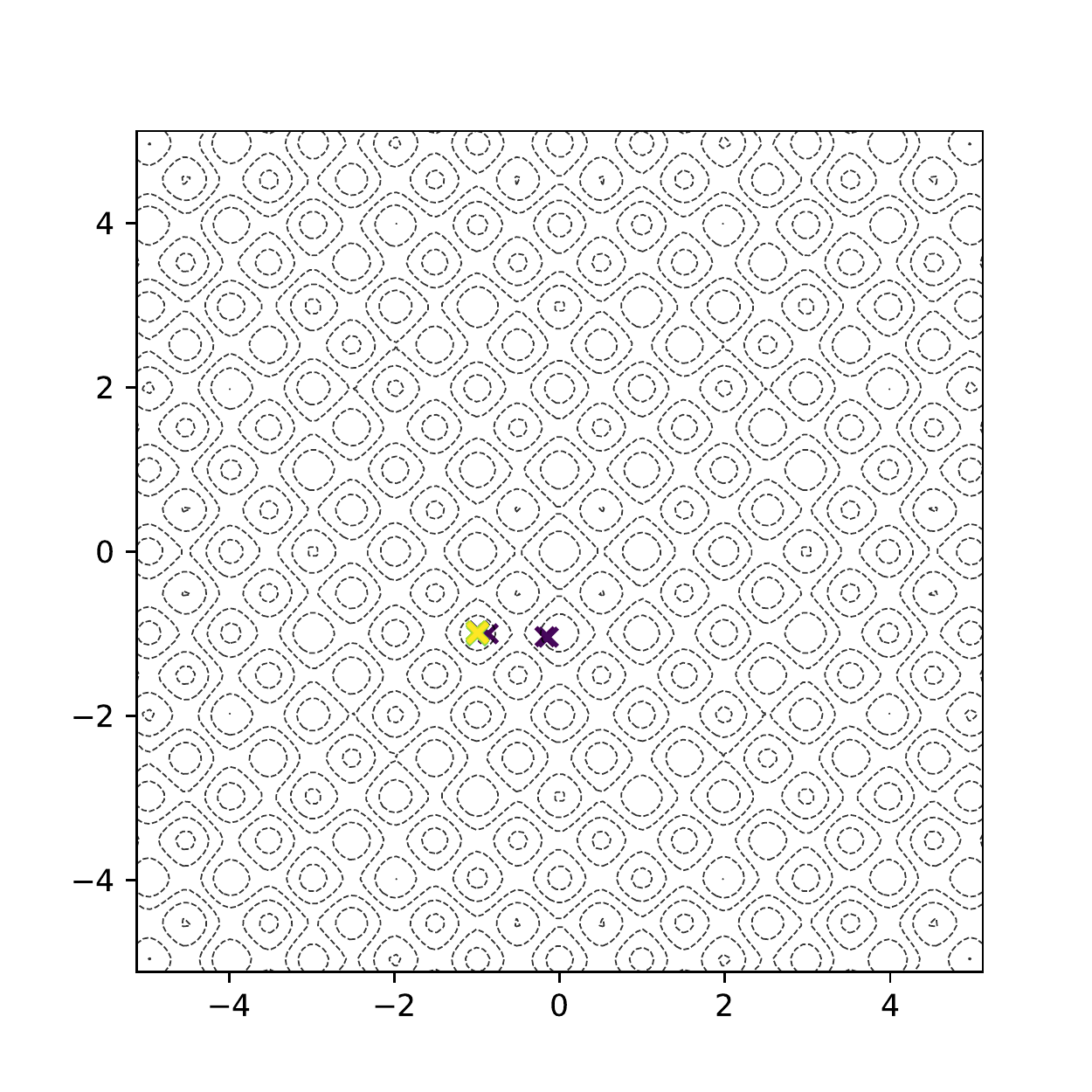}
\end{minipage}
\begin{minipage}[b]{0.29\textwidth}
\includegraphics[width=\textwidth,trim={0 .5cm 0 1.3cm},clip]{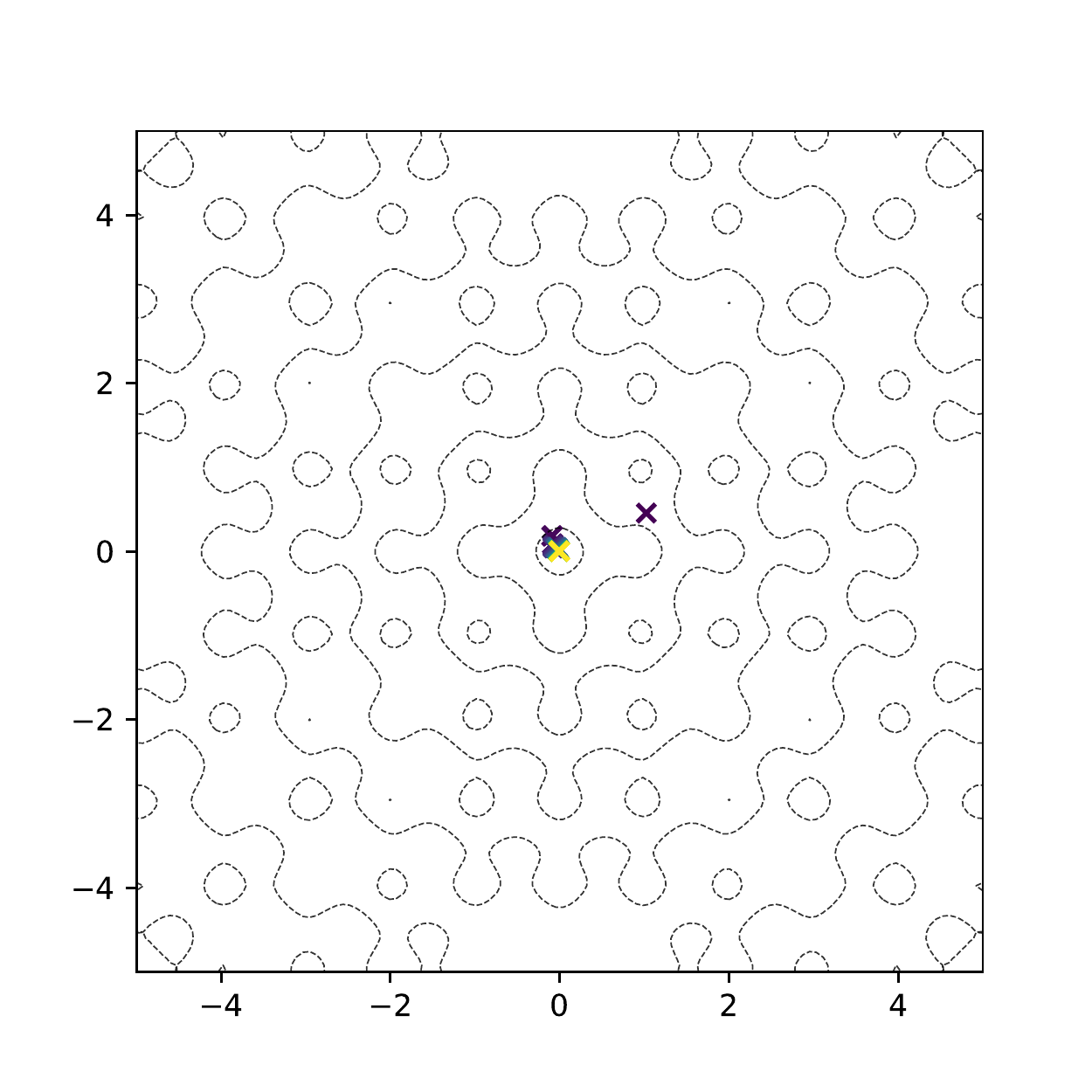}
\includegraphics[width=\textwidth,trim={0 .5cm 0 1.3cm},clip]{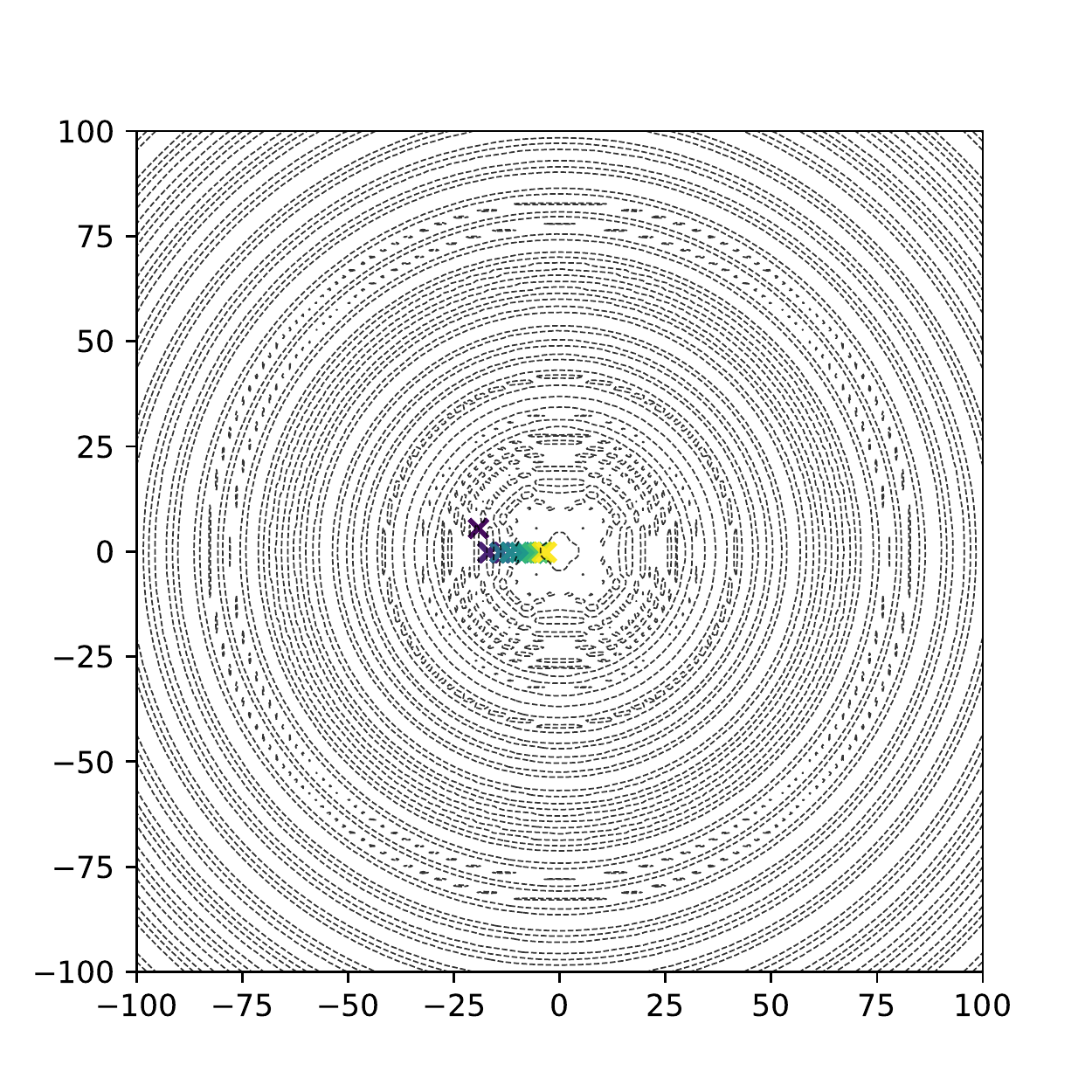}
\end{minipage}
\begin{minipage}[b]{0.29\textwidth}
\includegraphics[width=\textwidth,trim={0 .5cm 0 1.3cm},clip]{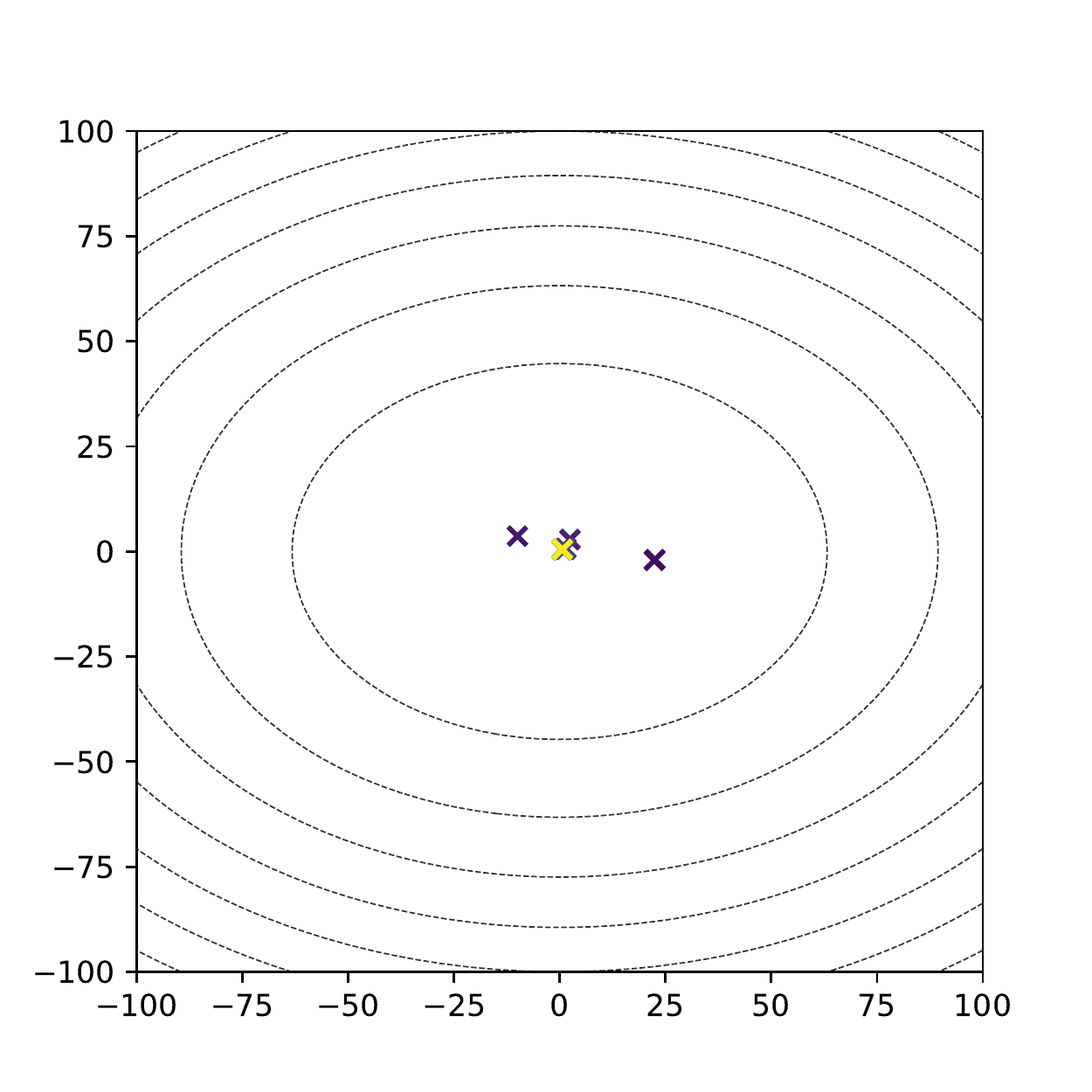}
\includegraphics[width=\textwidth,trim={0 .5cm 0 1.3cm},clip]{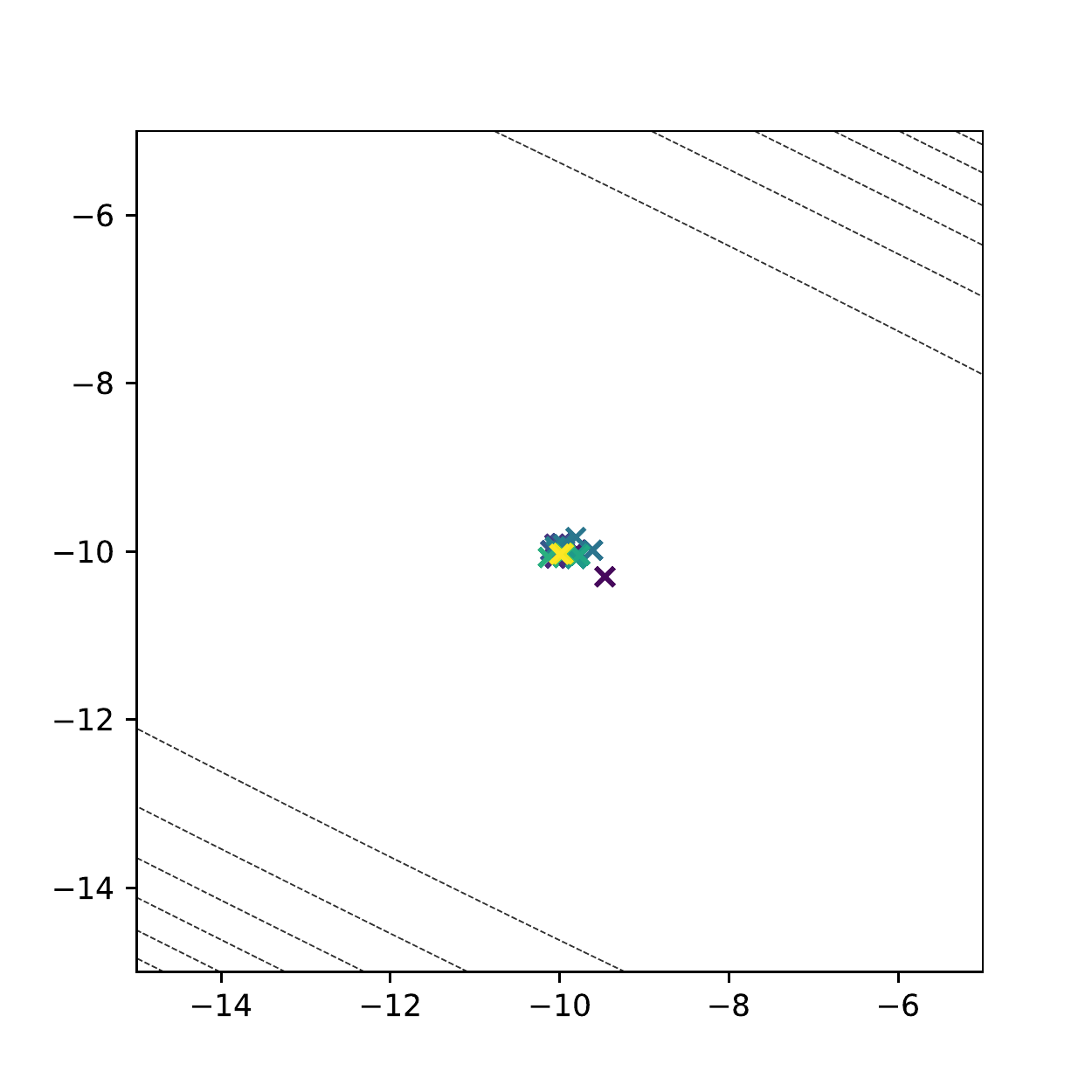}
\end{minipage}
\begin{minipage}[b]{0.08\textwidth}
\includegraphics[width=\textwidth,trim={0 .5cm 0 1cm},clip]{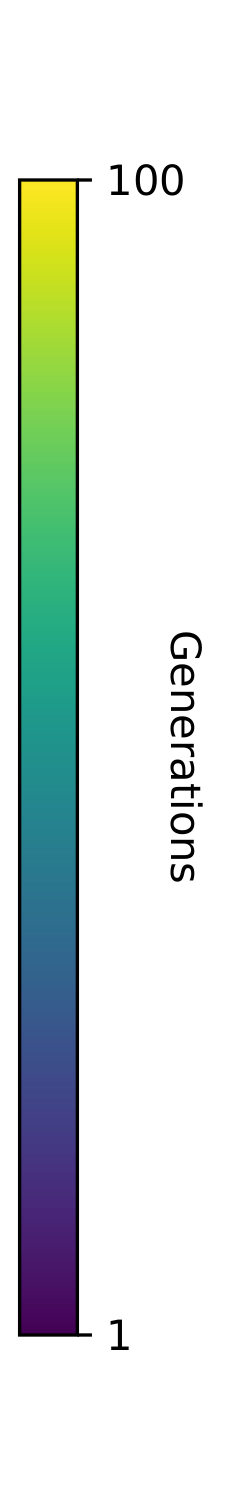}
\end{minipage}
\caption{\RA{Results of the single optimization run using GA for the $2$-dimensional
	benchmark functions. The colored
	crosses indicate the spatial coordinates of the best individual at each
	generation. Black lines indicate the isolines of the
	functions.\label{fig:conv_ga}}}
\end{figure}

\begin{figure}[h!]
\centering
\begin{minipage}[b]{0.29\textwidth}
\includegraphics[width=\textwidth,trim={0 .5cm 0 1.3cm},clip]{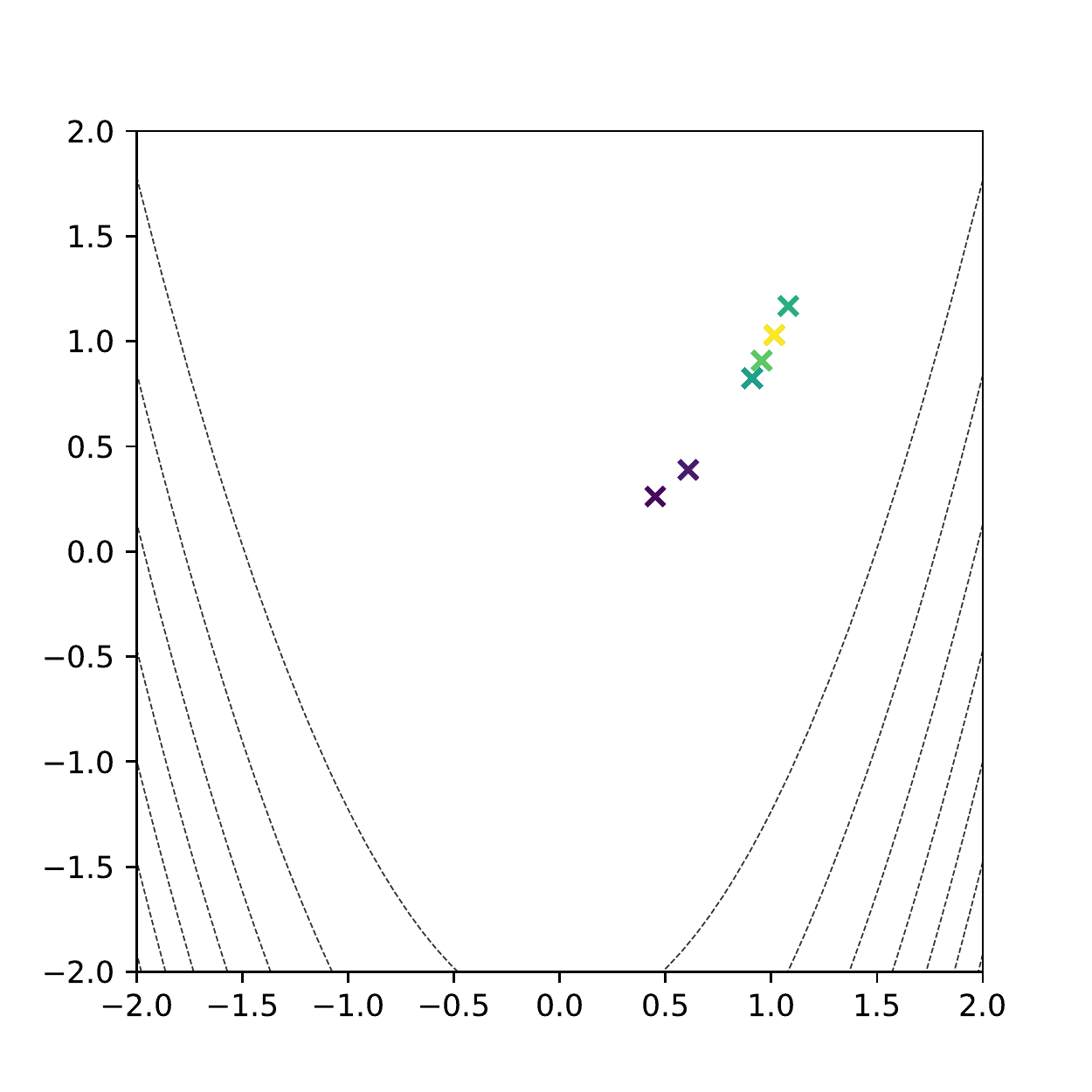}
\includegraphics[width=\textwidth,trim={0 .5cm 0 1.3cm},clip]{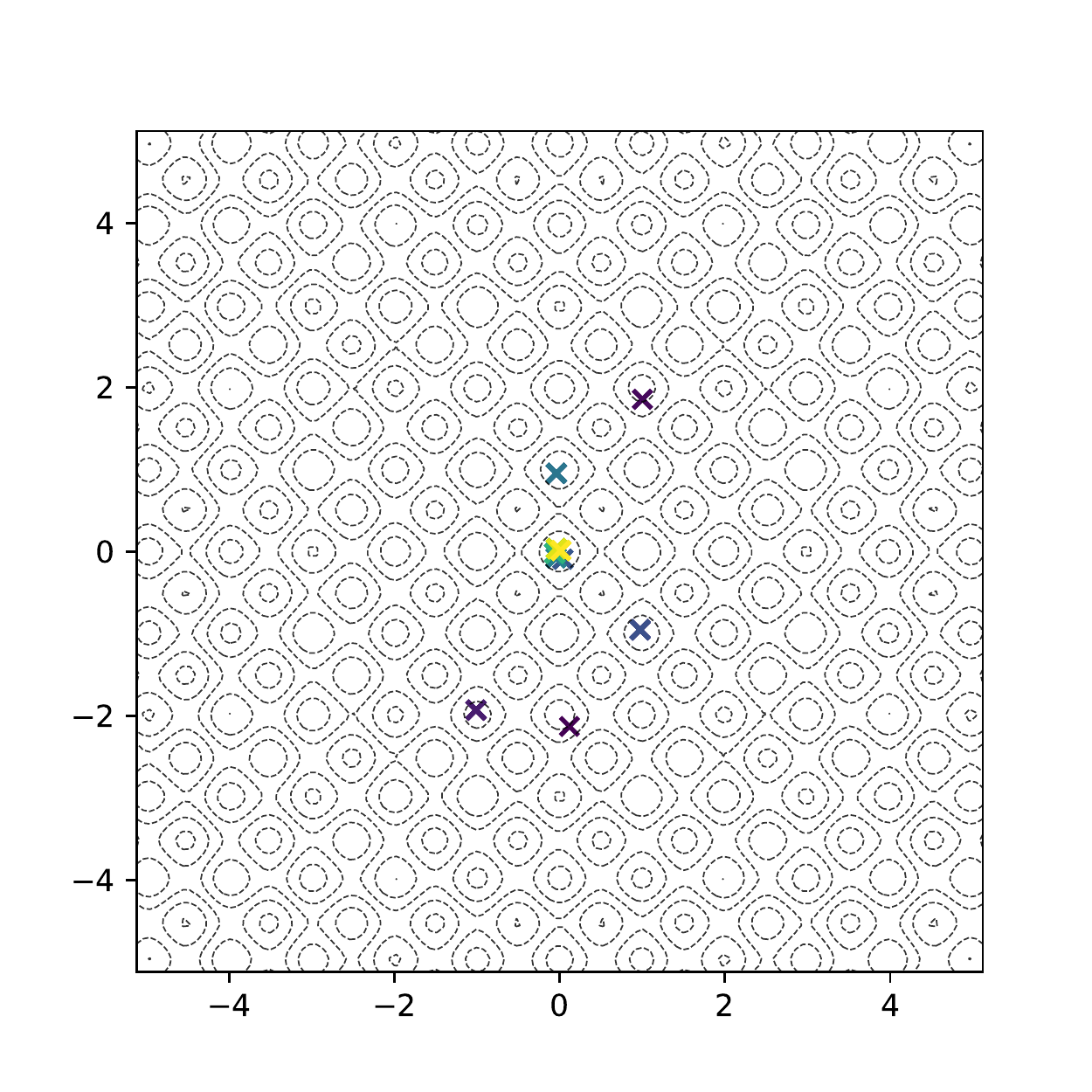}
\end{minipage}
\begin{minipage}[b]{0.29\textwidth}
\includegraphics[width=\textwidth,trim={0 .5cm 0 1.3cm},clip]{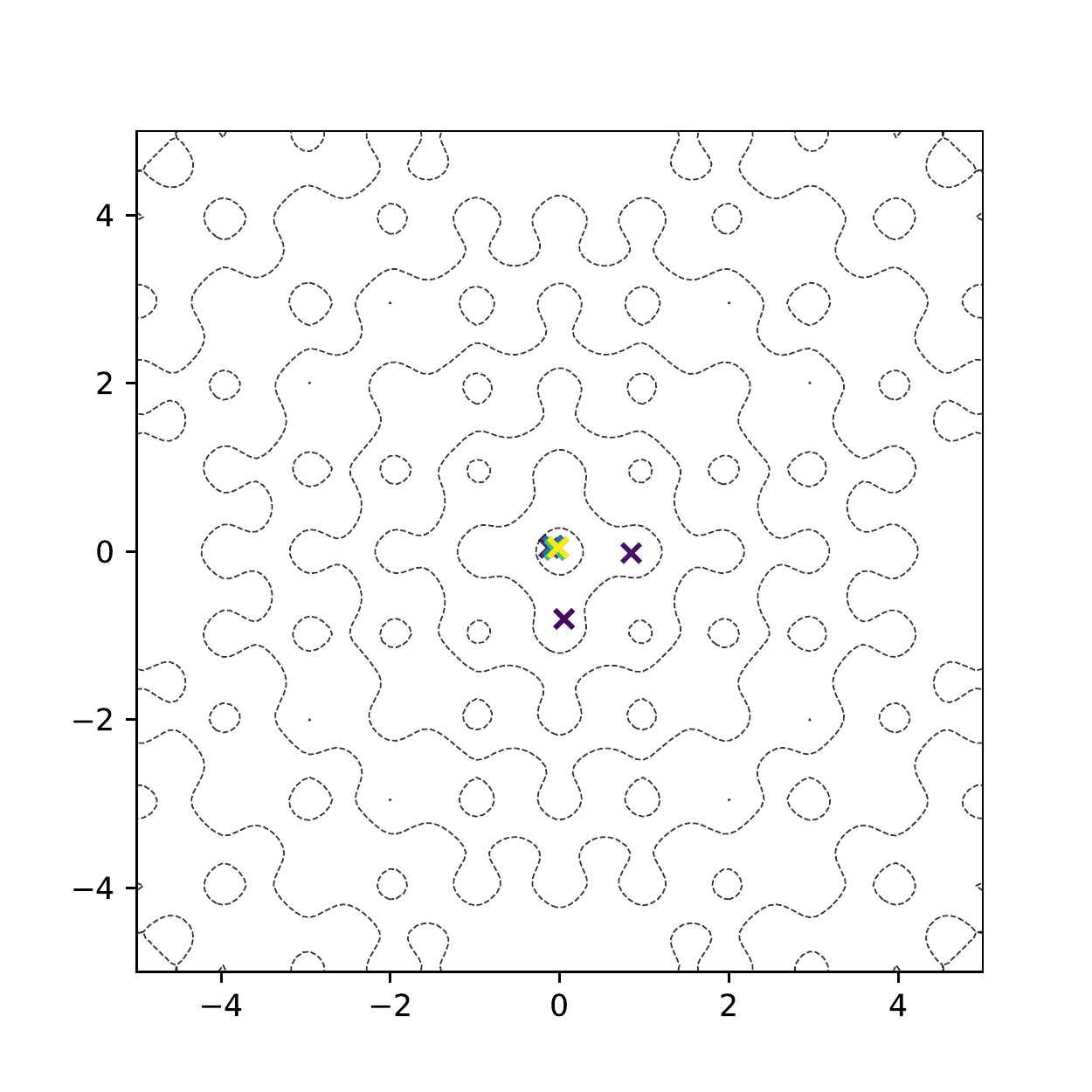}
\includegraphics[width=\textwidth,trim={0 .5cm 0 1.3cm},clip]{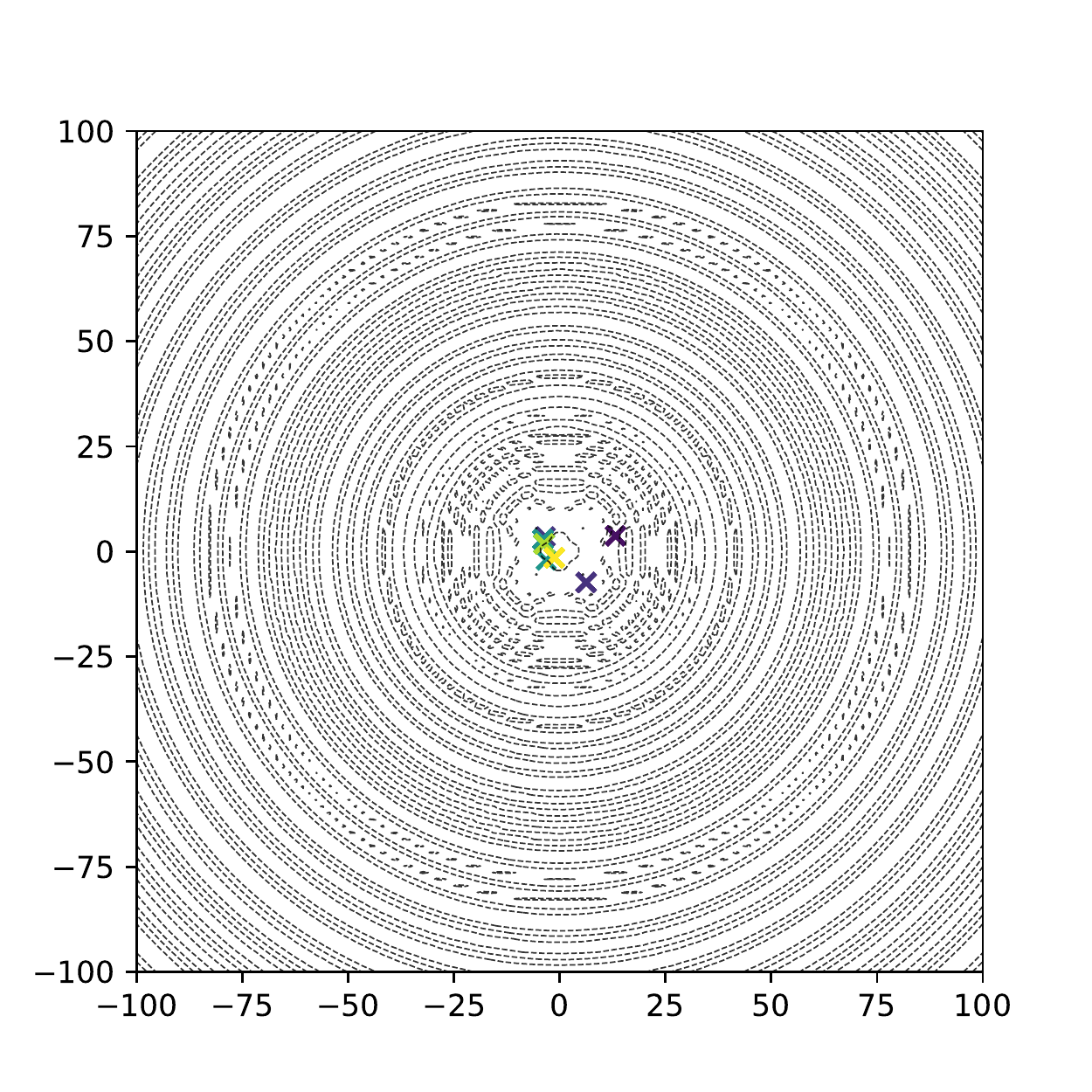}
\end{minipage}
\begin{minipage}[b]{0.29\textwidth}
\includegraphics[width=\textwidth,trim={0 .5cm 0 1.3cm},clip]{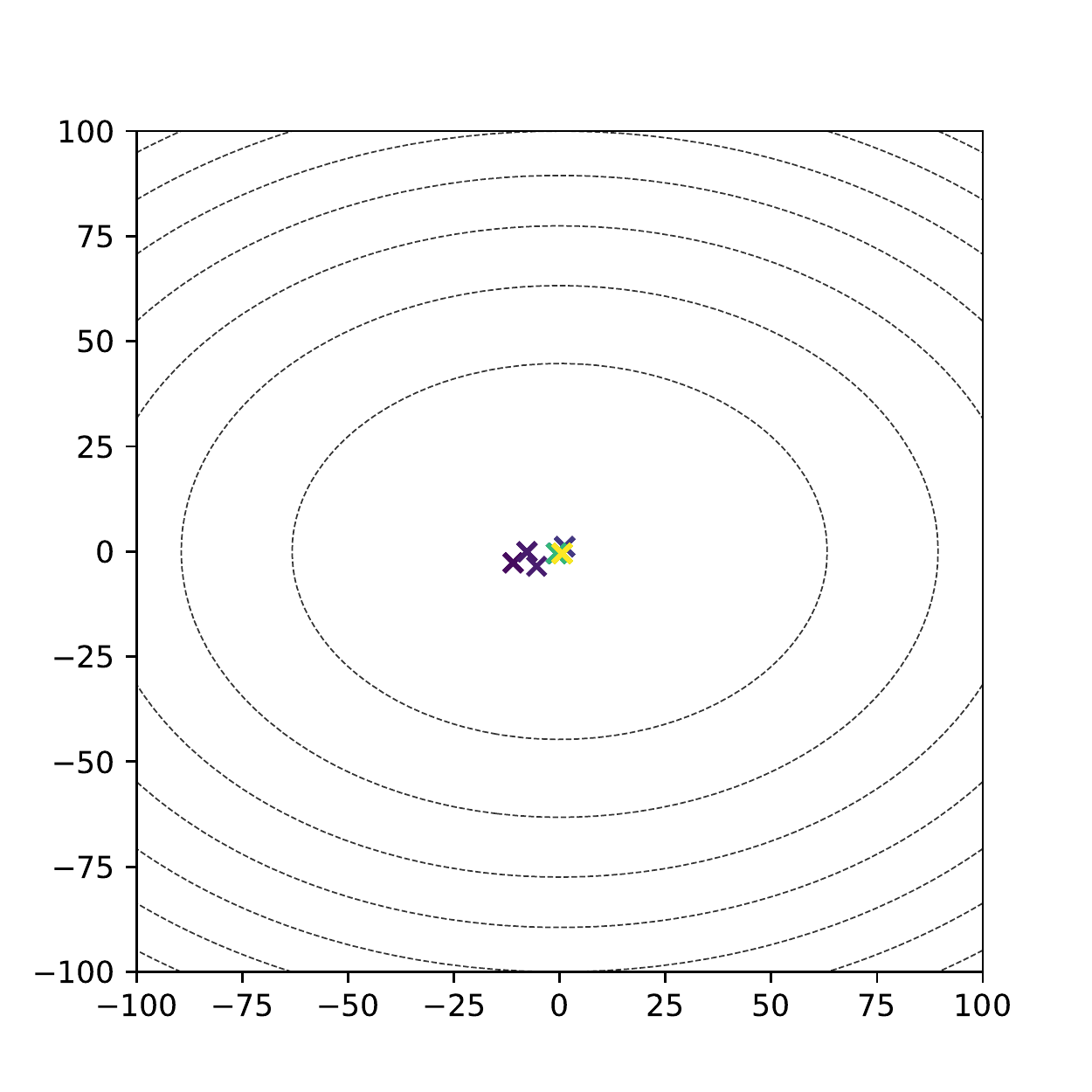}
\includegraphics[width=\textwidth,trim={0 .5cm 0 1.3cm},clip]{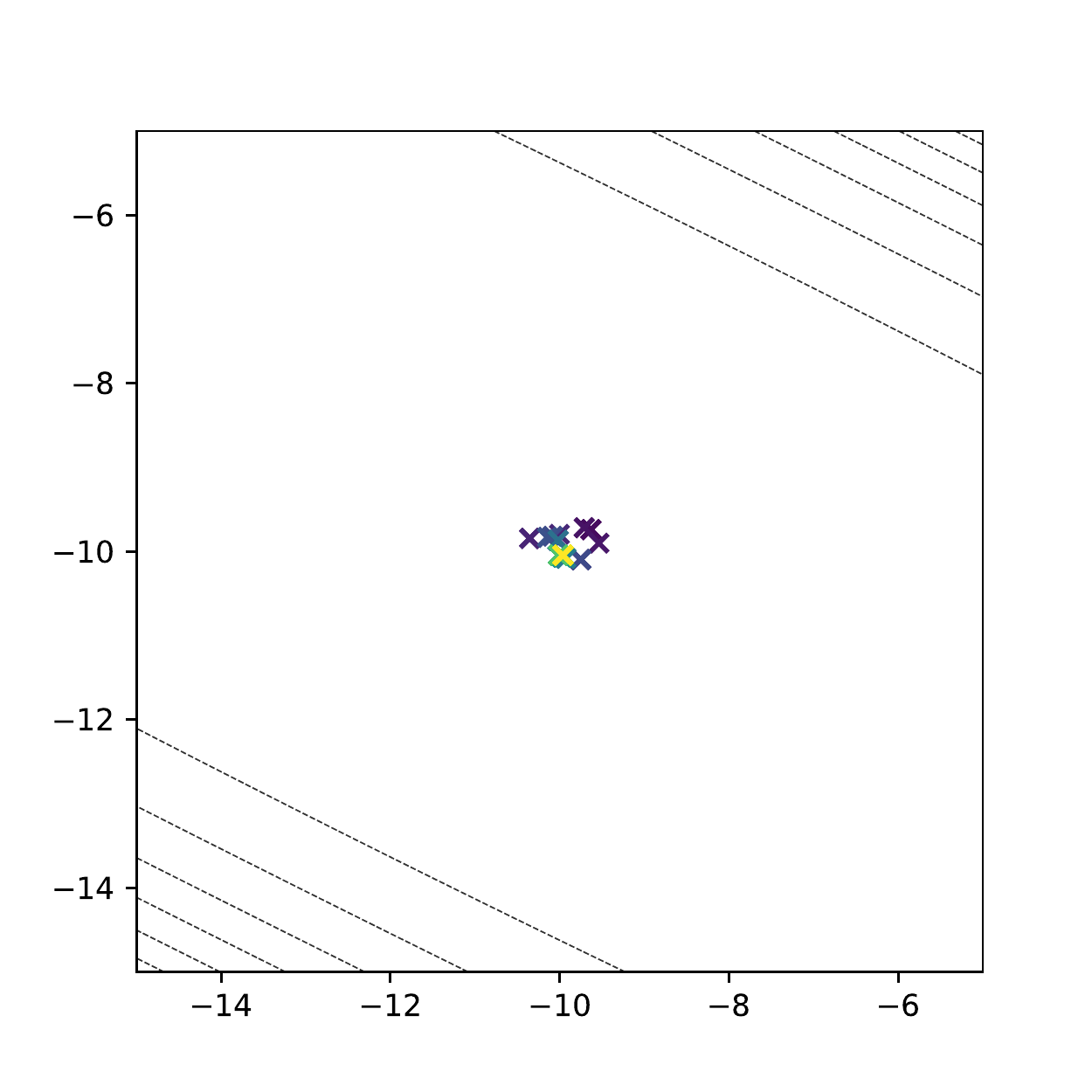}
\end{minipage}
\begin{minipage}[b]{0.08\textwidth}
\includegraphics[width=\textwidth,trim={0 .5cm 0 1cm},clip]{figures/colorbar.png}
\end{minipage}
\caption{\RA{Results of the single optimization run using ASGA
	for the $2$-dimensional benchmark functions. The colored
	crosses indicate the spatial coordinates of the best individual at each
	generation. Black lines indicate the isolines of the
	functions.\label{fig:conv_asga}}}
\end{figure}

We summarize in \cref{fig:conv_ga} and \cref{fig:conv_asga} the spatial
coordinates of the best individual after each generation using the standard GA
and ASGA. The proposed method reaches the global minimum for all the testcases,
performing better than the standard counterpart for the Rosenbrock and
Rastrigin functions. 

\begin{figure}[h!]
\centering
\includegraphics[width=.98\textwidth]{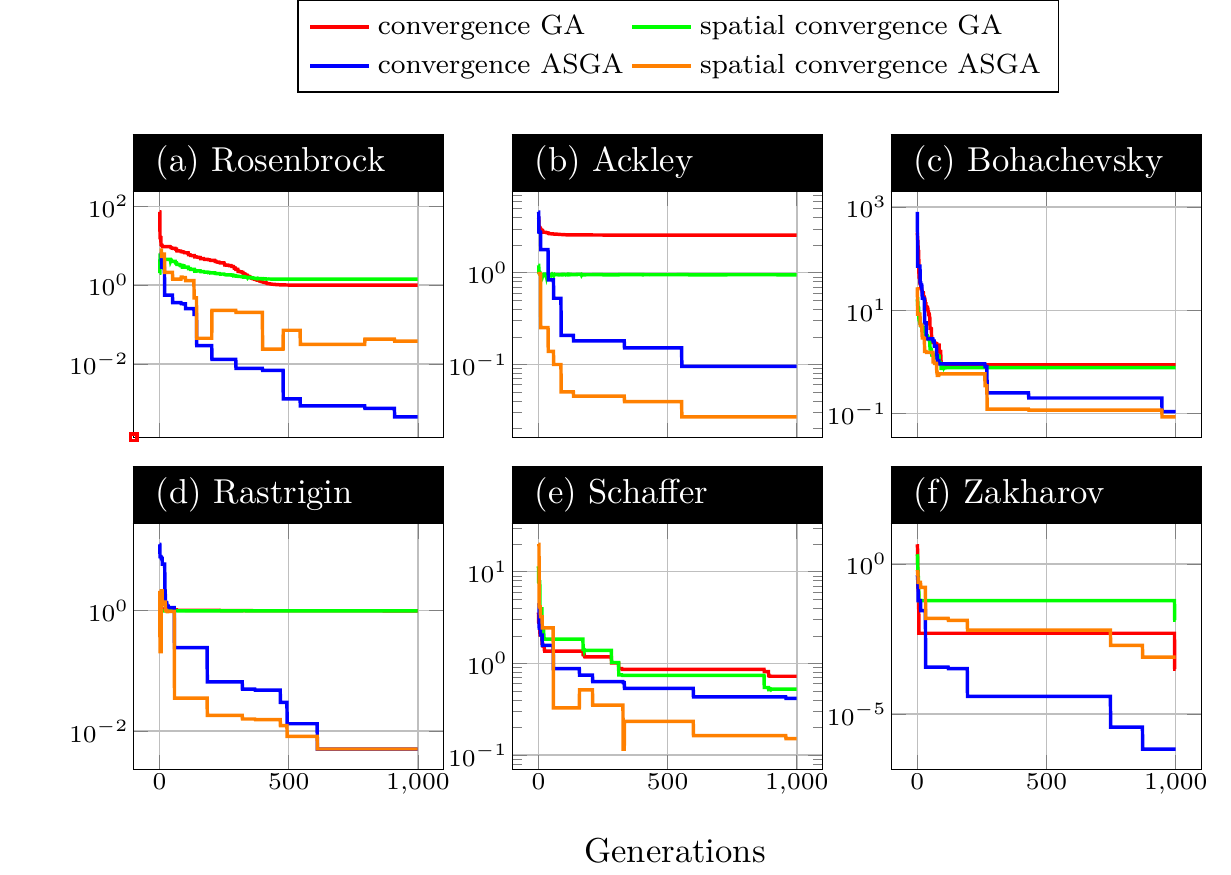}
	\caption{\RA{Convergence of GA and ASGA for the $2$-dimensional benchmark functions.\label{fig:conv}}}
\end{figure}

We also measure the convergence as the Euclidean distance between
the best individual fitness and the global optimum, and the \textit{spatial} convergence
as the Euclidean distance between the coordinates of the best individual and the
coordinates of the optimal point. We kept the same numerical settings, only
raising the number of generation to $1000$. \Cref{fig:conv} presents the
plots where we compare the trend using GA and ASGA: the proposed
method shows a better performance, not only thanks to the faster
convergence but also because in all the cases ASGA is able to get
closer than the GA to the global optimum.

\section*{Acknowledgements}
We thank Francesco Romor for the productive discussions and comments.
This work was partially supported by an industrial Ph.D. grant sponsored
by Fincantieri S.p.A. (IRONTH Project), by the MIUR 
through the FARE-X-AROMA-CFD project, by the INdAM-GNCS 2019 project ``Advanced
intrusive and non-intrusive model order reduction techniques and
applications'', and partially funded by European Union Funding for
Research and Innovation --- Horizon 2020 Program --- in the framework
of European Research Council Executive Agency: H2020 ERC CoG 2015
AROMA-CFD project 681447 ``Advanced Reduced Order Methods with
Applications in Computational Fluid Dynamics'' P.I. Professor Gianluigi Rozza.

\bibliographystyle{abbrv}

\end{document}